\documentclass{mcom-l}
\usepackage{amsmath}
\usepackage{amsthm}
\usepackage{epsfig}
\usepackage{hyperref}
\usepackage{amsfonts}
\usepackage{color}

\allowdisplaybreaks

\newtheorem{theorem}{Theorem}[section]
\newtheorem{lemma}[theorem]{Lemma}
\newtheorem{proposition}[theorem]{Proposition}

\theoremstyle{definition}

\theoremstyle{remark}

\numberwithin{equation}{section}

\newcommand{\nn}{\nonumber}
\def \endproof{\vrule height8pt width 5pt depth 0pt}
\def\refe#1{(\ref{#1})}

\def\d{\delta}

\def\R{\mathbb{R}}
\def\C{\mathbb{C}}

\def\d{{\rm d}}

\begin{document}

\title[Maximal regularity of FEMs for parabolic equations]{Maximal 
$\bf L^p$ analysis of finite element solutions for parabolic equations
with\\ nonsmooth coefficients in convex polyhedra}


\author{Buyang Li}
\address{{Department of Applied Mathematics, 
The Hong Kong Polytechnic University, Kowloon, Hong Kong.}\newline\indent 
Mathematisches Institut, 
Universit\"at T\"ubingen, 
D-72076 T\"ubingen, Germany. \newline\indent
Department of Mathematics,
Nanjing University, Nanjing 210093, P.R. China.
}  
\curraddr{}
\email{li@na.uni-tuebingen.de}

\thanks{The work of B. Li was supported in part by
NSFC (grant no. 11301262), 
and the research stay of the author at Universit\"at T\"ubingen
was supported by the Alexander von Humboldt Foundation}

\author{Weiwei Sun}
\address{Department of Mathematics,
City University of Hong Kong, 
Hong Kong.}
\curraddr{}
\email{maweiw@math.cityu.edu.hk}
\thanks{The work of W. Sun was supported in part by a grant
from the Research Grants Council of the Hong Kong SAR, China (project no. CityU 11301915)}
\thanks{This paper was accepted for publication in Math. Comp. 
on December 2, 2015.}

\subjclass[2010]{Primary 65M12, 65M60, Secondary 35K20} 

\date{}

\dedicatory{}

\begin{abstract}
The paper is concerned with Galerkin finite element solutions
of parabolic equations in a convex polygon or polyhedron with
a diffusion coefficient in $W^{1,N+\alpha}$ for some $\alpha>0$,
where $N$ denotes the dimension of the domain.
We prove the analyticity of the semigroup generated by the
discrete elliptic operator,
the discrete maximal $L^p$ regularity and the optimal $L^p$ error estimate
of the finite element solution for the parabolic equation.
\end{abstract}

\maketitle


\bibliographystyle{amsplain}

\section{Introduction}
\setcounter{equation}{0}

Let $\Omega$ be a bounded domain in $\R^N$
(with $N=2$ or $N=3$), and let $S_h$ be a finite element subspace of
$H^1_0(\Omega)$ consisting of continuous piecewise polynomials of
degree $r\geq 1$ subject to certain quasi-uniform triangulation of
the domain $\Omega$.
We consider the parabolic equation
\begin{align}\label{PDE0}
\left\{
\begin{array}{ll}
\partial_tu- \nabla \cdot (a \nabla u) =f
&\mbox{in}~~\Omega\times(0,\infty),\\
\displaystyle  u=0
&\mbox{on}~~\partial\Omega\times(0,\infty),\\
u(\cdot,0)=u^0 &\mbox{in}~~\Omega,
\end{array}
\right.
\end{align}
and its finite element approximation
\begin{align}\label{FEEq0}
\left\{
\begin{array}{ll}
\big(\partial_tu_h,v_h\big)
+ (a\nabla u_h, \nabla v_h) 
=(f,v_h) , ~~\forall~v_h\in S_h,
\\
u_h(0)=u^0_h ,
\end{array}
\right.
\end{align}
where $f$ is a given function,
and $a=(a_{ij}(x))_{N\times N}$ is an $N \times N$ symmetric matrix which
satisfies the ellipticity condition
\begin{align}\label{coeffcond}
\Lambda^{-1}|\xi|^2\leq
\sum_{i,j=1}^Na_{ij}(x)\xi_i\xi_j
\leq \Lambda|\xi|^2
,\quad\mbox{for}~~x\in\Omega,
\end{align}
for some positive constant $\Lambda$.

If we define the elliptic operator $A: H^1_0(\Omega)\rightarrow H^{-1}(\Omega)$
and its finite element approximation $A_h: S_h\rightarrow S_h$ by
\begin{align}
& (Aw,v):= (a \nabla w, \nabla v),
\qquad\qquad\,\, \forall~w,v\in
H^1_0(\Omega),
\label{A-e}
\\
 &(A_hw_h,v_h):= (a\nabla w_h, \nabla v_h) 
,\qquad \forall~w_h,v_h\in S_h ,
\label{Ah-e}
\end{align}
then the solutions of \refe{PDE0} and \refe{FEEq0} can be expressed by
\begin{align}
&u(t)=E (t)u^0 +\int_0^t E (t-s)f(s)\d s ,\\
&u_h(t)=E_h(t)u^0_h+\int_0^t E_h(t-s)f(s)\d s  ,
\end{align}
where $\{E(t)=e^{-tA}\}_{t>0}$ and
$\{E_h(t)=e^{-tA_h}\}_{t>0}$ denote the semigroups generated by the
operators $-A$ and $-A_h$, respectively.
By the theory of parabolic equations and \cite{Ouhabaz},
it is well known that
$\{E(t)\}_{t>0}$ is an analytic semigroup on
$C_0(\overline\Omega)$ satisfying
\begin{align}
&\|E(t)v\|_{L^\infty}
+t\|\partial_tE(t)v\|_{L^\infty}\leq C\|v\|_{L^\infty}
,\quad\forall\, v\in C_0(\overline\Omega), ~ \forall\, t>0 ,
\end{align}
which is equivalent to the resolvent estimate
\begin{align*}
&\|(\lambda+A)^{-1}v\|_{L^\infty} \leq C\lambda^{-1}\|v\|_{L^\infty}
,\quad\forall\, v\in C_0(\overline\Omega) ,
 \,\, \forall\, \lambda\in\Sigma_{\theta+\pi/2} ,
\end{align*}
where $\Sigma_{\theta+\pi/2}:=\{z\in \C : |{\rm arg}(z)|<\theta+\pi/2\}$.
The counterparts of these two inequalities above
for the discrete finite element
operator $A_h$ are the analyticity of the semigroup $\{E_h(t)\}_{t>0}$
on $L^\infty\cap S_h$:
\begin{align}
\|E_h(t)v_h\|_{L^\infty} +t\|\partial_tE_h(t)v_h\|_{L^\infty}\leq
C\|v_h\|_{L^\infty} ,\quad\forall\, v_h\in S_h,\, \forall\, t>0 ,
\label{STLEst}
\end{align}
and the resolvent estimate
\begin{align*}
\|(\lambda+A_h)^{-1}v_h\|_{L^\infty} \leq C\lambda^{-1}\|v_h\|_{L^\infty}
,\quad\forall\, v_h\in S_h , \,\, \forall\, \lambda\in\Sigma_{\varphi+\pi/2} .
\end{align*}
The estimates of the discrete semigroup
have attracted much attention in the past
several decades. With these estimates,
one may reach more precise analyses
of finite element solutions, such as
maximum-norm analysis of FEMs
\cite{Lucas,Tho, TW1,Wah},  error estimates of fully discrete FEMs
\cite{LN, Pal,Tho} and the discrete maximal $L^p$ regularity for
parabolic finite element equations \cite{Gei1,Gei2,KLL,Li,LS2}.

The proof of \refe{STLEst} dates back to Schatz et. al. \cite{STW1},
who proved \refe{STLEst} with a logarithmic factor
for the heat equation in
a two-dimensional smooth convex
domain with the linear finite element method.
The logarithmic factor was removed in the
case $r\geq 4$ for $N=1,2,3$ in \cite{NW},
and the analysis was further extended to the case
$1\leq N\leq 5$ in \cite{Chen}.
Later, a unified approach was presented in
\cite{STW98} by Schatz et. al., where they proved \refe{STLEst}
with the Neumann boundary condition for all
$r\geq 1$ and $N\geq 1$. The result was extended to
the Dirichlet boundary condition in \cite{TW2}
for the linear finite element method.
Some other maximum-norm error estimates can be found in
\cite{DLSW,DM,ELWZ,Han,Ley,Lin1}, and the resolvent estimates
can be found in \cite{Bak, BTW}.

A related topic is the discrete maximal $L^p$ regularity
(when $u^0=0$ and $1<p,q<\infty$)
\begin{align}
&\|\partial_tu_h \|_{L^p((0,T);L^q)}+\|A_hu_h \|_{L^p((0,T);L^q)}
\leq C_{p,q}\|f\|_{L^p((0,T);L^q)} ,
\label{LpqSt3}
\end{align}
which resembles the maximal $L^p$ regularity
of the continuous parabolic problem and was
proved by Geissert \cite{Gei1,Gei2}.
A straightforward
application of \refe{LpqSt3} is the $L^p$-norm error estimate
\begin{align}
\| P_h u - u_h \|_{L^p((0,T);L^q)}
\le C_{p,q}(\| P_h u^0 - u_h^0 \|_{L^q}
+ \| P_h u - R_h u \|_{L^p((0,T);L^q)} )  ,
\label{error}
\end{align}
where $R_h$ is the Ritz projection
associated with the operator $A$
and $P_h$ is the $L^2$ projection
onto the finite element space.

All these estimates
were established under the assumption
that the coefficients $a_{ij}$ and the domain
$\Omega$ are smooth enough
so that the parabolic Green's function satisfies
\begin{align}\label{g-cond000}
|\partial_t^{\gamma}\partial_x^\beta G(t, x, y)|\leq
C(t^{1/2}+|x-y|)^{-(N+2\gamma+|\beta|)} e^{-\frac{|x-y|^2}{Ct}}
 ,~~\forall~0\leq \gamma\leq 2,\, 0\leq |\beta|\leq 2 \, .
\end{align}
Although the condition on the coefficients was relaxed to
$a_{ij} \in C^{2+\alpha}(\overline\Omega)$
in \cite{Gei1}, this assumption
is still too strong for many physical applications.
One of the examples is an
incompressible miscible flow in porous media \cite{Dou,LS1},
where the diffusion-dispersion tensor $[a_{ij}]_{i,j=1}^N$
is only a Lipschitz continuous function of the velocity field.
In a recent work \cite{Li}, the first author proved
\refe{STLEst} in a smooth domain
under the assumption $a_{ij} \in W^{1,\infty}(\Omega)$,
together with the estimate (when $u^0=0$
and $1<p,q<\infty$)
\begin{align}
\|u_h \|_{L^p((0,T);W^{1,q})}
\leq C_{p,q}\|f\|_{L^p((0,T);W^{-1,q})} ,
\label{LpqSt2}
\end{align}
which were then applied to the incompressible
miscible flow in porous media \cite{LS2}.
Moreover, the problem in a polygon or a
polyhedron is of high interest
in practical cases, while the
inequality \refe{g-cond000}
does not hold in arbitrary convex polygons or polyhedra,
and all the analyses of
{\rm\refe{LpqSt3}}-{\rm(\ref{LpqSt2})} are limited
to smooth domains so far.
For the problem in two-dimensional
polygons with constant coefficients,
the inequality (\ref{STLEst}) with an extra logarithmic factor
was proved in \cite{CLTW,Ran,Tho}
by using the following
estimate of the discrete Green's function $\Gamma_h$:
$$
\int_{\Omega} | \Gamma_h(t,x,x_0) | dx \le C |\ln h|
\, .
$$
The corresponding results in three-dimensional polyhedra are unknown.
More interested is whether these stability estimates hold with the natural regularity $a_{ij} \in W^{1,p}(\Omega)$ for some $1<p<\infty$, since
such estimates are important for the extension of the analysis
to a general nonlinear model.

This paper focuses on \refe{STLEst}-\refe{LpqSt3} and
\refe{LpqSt2} in a convex polygon or polyhedron
with a weaker regularity of the diffusion coefficient.
Instead of estimating $\Gamma_h$ directly,
we present a more precise estimate for the error function
$F:=\Gamma_h-\Gamma$ (see Lemma 2.2) with which the
logarithmic factor can be removed (this idea was used in \cite{STW98}),
where $\Gamma$ is a regularized Green's function.
To compensate the lack of pointwise estimate of
the second-order derivatives of the Green's function,
we use local $W^{1,\infty}$ estimate
and local energy estimates of the
second-order derivatives (see Lemma 4.1).
Our main result is the following theorem.
\begin{theorem}\label{MainTHM1}
{\it Assume that $a_{ij}\in W^{1,N+\alpha}(\Omega)$ for some $\alpha>0$,
satisfying the condition {\rm (\ref{coeffcond})}, and assume that
$\Omega$ is either a
convex polygon in $\R^2$ or a convex polyhedron in $\R^3$.
Then

(1) the semigroup estimate {\rm\refe{STLEst}} holds,

(2) the solution of {\rm(\ref{FEEq0})}
satisfies {\rm(\ref{LpqSt3})}
when $f\in L^p((0,T);L^q)$ and $u^0=0$, 

(3) the solution of {\rm(\ref{FEEq0})}
satisfies {\rm(\ref{LpqSt2})}
when $f\in L^p((0,T);W^{-1,q})$ and $u^0=0$.
}
\end{theorem}

Under the assumptions in Theorem \ref{MainTHM1}
and assuming that the solution of \refe{PDE0} satisfies
$u \in C(\overline\Omega\times[0,T])$,
\refe{error} follows immediately from {\rm(\ref{LpqSt3})}.

The rest of this paper is organized as follows. In section 2,
we introduce some notations and present a key lemma based on which
our main theorem can be proved.  In section 3, we present superapproximation
results for smoothly truncated finite element functions and present several
estimates for the parabolic Green's functions under the assumed
regularity of the coefficients and the domain.  Based on these estimates,
we prove our key lemma in section 4.

\section{Notations, assumptions and sketch of the proof}
\setcounter{equation}{0}

\subsection{Notations}
For any nonnegative integer $k $ and $1\leq p\leq\infty$,
we let $W^{k,p}(\Omega)$ be the conventional Sobolev space of functions
defined in $\Omega$, and let $W^{1,p}_0(\Omega)$ be the subspace of
$W^{1,p}(\Omega)$ consisting of functions whose traces vanish on $\partial\Omega$.
As conventions, we denote the dual space of $W^{1,p}_0(\Omega)$
by $W^{-1,p'}(\Omega)$ for $1\leq p<\infty$, and denote
$H^k(\Omega):=W^{k,2}(\Omega)$ and $L^p(\Omega):=W^{0,p}(\Omega)$ for
any integer $k$ and $1\leq p\leq \infty$.

Let $ Q_T:=\Omega\times(0,T)$.
For any Banach space $X$ and a given $T>0$,
we let $L^p((0,T);X)$ be the Bochner spaces equipped with the norm
\begin{align*}
&\|f\|_{L^p((0,T);X)} =\left\{
\begin{array}{ll}
\displaystyle\biggl(\int_0^T
\|f(t)\|_X^pdt\biggl)^\frac{1}{p},
&
1\leq p<\infty
,\\[10pt]
\displaystyle{\rm ess\,\,}\sup_{\!\!\!\!\!\!\!\!\!\! t\in(0,T)}\, \|f(t)\|_X .
& p=\infty,\end{array}
\right.
\end{align*}
To simplify notations, in the following sections, we write $L^p$, $H^k$
and $W^{k,p}$ as the abbreviations of $L^p(\Omega)$, $H^k(\Omega)$
and $W^{k,p}(\Omega)$, respectively, and denote by $(\cdot\, , \cdot\, )$ the inner product in $L^2(\Omega)$.
For any subdomain $Q\subset Q_T$, we define
\begin{align*}
&Q^t:=\{x\in\Omega:~ (x,t)\in Q\},\\[5pt]
&\|f\|_{L^{\infty,2}(Q)}:={\rm ess}\!\!\sup_{t\in(0,T)}
\|f(\cdot,t)\|_{L^2(Q^t)} ,\\
&
\|f\|_{L^p(Q)}:=\bigg(\iint_Q |f(x,t)|^p\d x\d t\bigg)^{\frac{1}{p}} ,
\quad\forall\, 1\leq p<\infty,
\end{align*}
and denote $w(t)=w(\cdot,t)$ for any function $w$ defined on $ Q_T$.

We assume that $\Omega$ is partitioned into quasi-uniform triangular
elements $\tau_l^h$, $l=1,\cdots,L$, with  $h=\max_{l}\{\mbox{diam}\,\tau_l^h\}$,
and let $S_h$ be a finite element subspace of $H^1_0(\Omega)$ consisting
of continuous piecewise polynomials of degree $r\geq 1$ subject to
the triangulation.
Let $a(x)=\left(a_{ij}(x) \right)_{N\times N}$ be the coefficient matrix
and define the operators
\begin{align*}
&A:H^1_0\rightarrow H^{-1},\qquad~\, A_h:S_h\rightarrow S_h,\\
&R_h:H^1_0\rightarrow S_h, \qquad~~\, P_h:L^2\rightarrow S_h ,
\end{align*}
by
\begin{align*}
&\big(A\phi,v\big)= \big(a\nabla \phi,\nabla
v\big) &&
\mbox{for all~~$\phi ,v \in H^1_0$},\\[3pt]
&\big(A_h\phi_h,v\big)= \big(a\nabla \phi_h,\nabla
v\big) &&
\mbox{for all~~$\phi_h\in S_h$, $v\in S_h$},\\[3pt]
&\big(A_hR_hw,v\big)=\big(Aw,v\big)  &&\mbox{for
all~~$w \in H^1_0$ and $v\in S_h$} ,
\\[3pt]
&\big(P_h\phi,v\big)=\big(\phi,v\big) &&
\mbox{for all~~$\phi \in L^2$
and $v \in S_h$} .
\end{align*}
Clearly, $R_h$ is the Ritz projection operator associated to the elliptic
operator $A$ and $P_h$ is the $L^2$ projection operator onto the finite element
space.
The following estimates are useful in this paper.

\begin{lemma}\label{LemMaxLp}
{\it
If $\Omega$ is a bounded convex domain
and $a_{ij}\in W^{1,N+\alpha}(\Omega)$, $N\geq 2$,
then we have
\begin{align}
&\|w\|_{H^2}\leq
C\| \nabla\cdot(a\nabla w)\|_{L^2},
&&\forall\, w\in H^1_0,
\label{H2Reg0}\\
&\|\nabla w\|_{L^\infty}
\leq C_{p}\| \nabla\cdot(a\nabla w)\|_{L^p}, \qquad
\mbox{for any given $p>N$},
&&\forall\, w\in H^1_0 ,
\label{W1inftyReg}
\end{align}
and the  solution of \refe{PDE0}
with $u^0=0$ satisfies
\begin{align}
&\|\partial_tu\|_{L^p((0,T);L^q)}
+\|Au\|_{L^p((0,T);L^q)}
\leq C_{p,q}\|f \|_{L^p((0,T);L^q)} ,
\label{LpW2quf}\\
&\|\partial_tu\|_{L^p((0,T);W^{-1,q})}
+\|u\|_{L^p((0,T);W^{1,q})}
\leq C_{p,q}\|f \|_{L^p((0,T);W^{-1,q})},
\label{LpW1quf}
\end{align}
for all $1<p,q<\infty$. 
}
\end{lemma}

In the Lemma above, \refe{H2Reg0} is the
standard $H^2$-regularity estimate in convex domains and
\refe{W1inftyReg} is a simple consequence of
the Green's function estimates given
in Theorem 3.3--3.4 of \cite{GW}, 
and \refe{LpW2quf}-\refe{LpW1quf} are consequences of the
maximal $L^p$ regularity (see Appendix for details).

\subsection{Properties of the finite element space and Green's functions}
\label{Sec2-2}
For any subdomain $D\subset\Omega$,
we denote by $S_h(D)$ the space of functions restricted to
the domain $D$, and denote by $S_h^0(D)$ the subspace of $S_h(D)$ consisting of functions which equal zero
outside $D$. For any given
subset $D\subset\Omega$, we denote
$B_d (D) =\{x\in\Omega: {\rm dist}(x,D)\leq d\}$ for $d>0$.
Then there exist positive constants
$K $ and $\kappa$ such that the
triangulation and the corresponding
finite element space $S_h$ possess the following properties
($K$ and $\kappa$ are independent of the subset $D$ and $h$).
\medskip

{\bf(P0)$\,$ Quasi-uniformity:}

For all triangles (or tetrahedron) $\tau_l^h$ in the partition,
the diameter $h_l$ of $\tau_l^h$ and the radius $\rho_l$
of its inscribed ball satisfy
$$
K^{-1}h\leq \rho_l \leq h_l\leq Kh .
$$

{\bf(P1)$\,$ Inverse inequality:}

If $D$ is a union of elements in the partition, then
\begin{align*}
\|\chi_h\|_{W^{l,p}(D)}
\leq K h^{-(l-k)-(N/q-N/p)}\|\chi_h\|_{W^{k,q}(D)} ,
\quad\forall\,\,\chi_h\in S_h,
\end{align*}
for $0\leq k\leq l\leq 1$ and 
$1\leq q\leq p\leq\infty$.  

{\bf(P2)$\,$ Local approximation and superapproximation:}

(1) There exists a linear operator
$I_h:H^1_0(\Omega)\rightarrow S_h$
such that if $d\geq \kappa h$, then
\begin{align*}
&\|v-I_hv\|_{L^2(D)} 
\leq K\sum_{l=0}^kh^{k}d^{-l}\|v\|_{H^{{k-l}}(B_d(D))} ,
\quad \forall\,\,v\in H^k\cap H^1_0 , \,\,\,
1\leq k\leq 2 .
\end{align*}
Moreover, if supp$(v)\subset \overline D$,
then $I_hv\in S_h^0(B_{d}(D))$.
For example, the Cl\'ement interpolation operator
defined in \cite{Clement} has these properties.
Also, the Lagrange interpolation operator $\Pi_h$ satisfies
\begin{align*}
&\|v-\Pi_hv\|_{L^2(D)}
+h\|\nabla(v-\Pi_hv)\|_{L^2(D)}
\leq Kh^2\|\nabla^2v\|_{L^2(B_d(D))} ,
\quad \forall\,\,v\in H^2\cap H^1_0  .
\end{align*}

(2) If $d\geq \kappa h$,
$\omega=0$ outside $B_{2d}(D)$
and $|\partial^\beta\omega|\leq Cd^{-|\beta|}$
for all multi-index $\beta$,
then for any $\psi_h\in S_h(B_{3d}(D))$
there exists $\eta_h\in S_h^0(B_{3d}(D))$ such that
\begin{align*}
&\|\omega\psi_h-\eta_h\|_{H^k(B_{3d}(D))}
\leq K  h^{1-k}d^{-1}
\|\psi_h\|_{L^2(B_{3d}(D))} ,
\quad k=0,1 .
\end{align*}
Furthermore, if $\omega\equiv 1$
on $B_{d}(D)$, then $\eta_h=\psi_h$ on $D$ and
\begin{align*}
&\|\omega\psi_h-\eta_h\|_{H^k(B_{3d}(D))}
\leq K h^{1-k}d^{-1}
\|\psi_h\|_{L^2(B_{3d}(D)\backslash D)} ,
\quad k=0,1 .
\end{align*}
For example, $\eta_h=\Pi_h(\omega\psi_h)$ has
these properties.\medskip

{\bf(P3)$\,$ Regularized Delta function:}

For any $x_0\in\overline\tau_j^h$, there exists a function
$\widetilde\delta_{x_0}\in C^3(\overline\Omega)$ with support in $\tau_j^h$
such that
\begin{align*}
&\chi_h(x_0)=\int_{\tau_j^h}\chi_h \widetilde\delta_{x_0}\d x,
\quad\forall\,\chi_h\in S_h ,\\
&\|\widetilde\delta_{x_0}\|_{W^{l,p}}
\leq K h^{-l-N(1-1/p)}
\quad\mbox{for}\,\,\,1\leq p\leq\infty,
\,\,\, l=0,1,2,3 .
\end{align*}

{\bf(P4)$\,$ Discrete Delta function}

Let $\delta_{x_0}$ denote the Dirac Delta function centered at
$x_0$, i.e. $\int_\Omega\delta_{x_0}(y)\varphi (y)\d y=\varphi(x_0)$ for
any $\varphi\in C(\overline\Omega)$.  The discrete Delta function
$P_h  \widetilde \delta_{x_0}$ satisfies that
\begin{align*}
& P_h  \widetilde\delta_{x_0}(x)
\leq Kh^{-N} e^{-\frac{|x-x_0| }{K h }} ,
\quad \forall\, x,x_0\in\Omega  .
\end{align*}

The properties (P0)-(P4) hold for any quasi-uniform partition with
those standard finite element spaces and also, have been used in
many previous works such as \cite{Li,STW98,SW95,TW2}.
The proof can be found in the appendix of \cite{SW95}.\medskip

For an element $\tau_l^h$ and a point $x_0\in \overline\tau_l^h$, we
let $G(t,x,x_0)$ be the Green's function of the parabolic
equation, defined by
\begin{align}\label{GFdef}
&\partial_tG(t,\cdot,x_0)+ AG(t,\cdot,x_0)=0\quad
\mbox{for $t>0$ with $G(0,x,x_0)=\delta_{x_0}(x)$} ,
\end{align}
The regularized Green's function $\Gamma(t,x,x_0)$ is defined by
\begin{align}\label{GMFdef}
&\partial_t\Gamma(\cdot,\cdot,x_0)+A\Gamma(\cdot,\cdot,x_0)=0\quad
\mbox{for~ $t>0$~ with~
$\Gamma(0,\cdot,x_0)=\widetilde\delta_{x_0}$},
\end{align}
where $\widetilde\delta_{x_0}$ is given in (P2), and the discrete
Green's function $\Gamma_h(\cdot,\cdot,x_0)$ is defined by
\begin{align}
&\partial_t\Gamma_{h}(\cdot,\cdot,x_0)
+ A_h\Gamma_{h}(\cdot,\cdot,x_0)=0 \label{GMhFdef}
\quad\mbox{for~ $t>0$~ with~
$\Gamma_{h}(0,\cdot,x_0)=P_h\widetilde\delta_{x_0}$} .
\end{align}
The functions $G(t,x,x_0)$ and $\Gamma_h(t,x,x_0)$
are symmetric with respect
to $x$ and $x_0$.

By the fundamental estimates of parabolic equations,
there exists a positive constant $C$ such that
(\cite{FS}, Theorem 1.6; note that the Green's function
in the domain $\Omega$ is less than the Green's function in $\R^N$)
\begin{align}
&|G(t,x,y)|\leq
C(t^{1/2}+|x-y|)^{-N}e^{-\frac{|x-y|^2}{Ct}} .
\label{FEstP}
\end{align}
By estimating $\Gamma(t,x,x_0)=\int_\Omega G(t,x,y)\widetilde\delta_{x_0}(y)\d y$, it is easy to see that (\ref{FEstP}) also holds
when $G$ is replaced by $\Gamma$ and
when $\max(t^{1/2},|x-y|)\geq 2h$.

\subsection{Decomposition of the domain $\Omega\times(0,T)$}
\label{SecGF}
Here we present some further notations
on a dyadic decomposition of the domain $\Omega\times(0,T)$, which were introduced in \cite{STW98} and also used in many other articles \cite{Gei1,Ley,Li,TW2}.
Let $R_0$ be the smallest distance between
a corner and a closed face which does not contained this corner.

For the given polygon/polyhedron $\Omega$,
there exists a positive constant $K_0\geq \max(1,R_0)$
(which depends on the interior angle of the edges/corners
of $\Omega$) such that

(1) if $z_0$ is a point in the interior of $\Omega$ and
$B_\rho(z_0)$
intersects a face of $\Omega$, then
$B_\rho(z_0)\subset B_{2\rho }(z_1)$ for
some $z_1$ which is on a face of $\Omega$;

(2) if $z_1$ is on a face of $\Omega$ and
$B_\rho(z_1)$
intersects another face, then
$B_\rho(z_1)\subset B_{\rho K_0}(z_2)$ for
some $z_2$ which is on an edge of $\Omega$;

(3) if $z_2$ is on an edge of $\Omega$ and
$B_\rho(z_2)$
intersects another face which does not contain this edge,
then $B_\rho(z_2)\subset B_{\rho K_0}(z_3)$ for
some $z_3$ which is a corner of $\Omega$.

\begin{figure}[ht]
\vspace{0.1in}
\centering
\epsfig{file=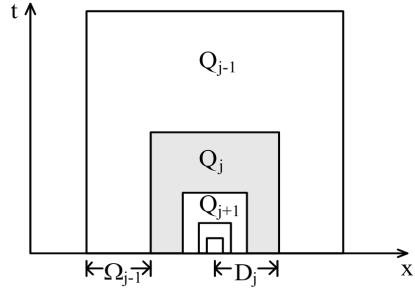,height=1.6in,width=2.1in}
\caption{Illustration of the
subdomains $Q_j$, $\Omega_j$ and $D_j$.}
\label{DomainFig}
\end{figure}

For any integer $j\geq 1$, we define $d_j=2^{-j-3}R_0K_0^{-2}$.
For a given $x_0\in\Omega$,
we let $J_*$ be an integer satisfying
$d_{J_*}=2^{-J_*-3}R_0K_0^{-2}= C_*h$
with $C_*\geq \max(10,10\kappa,R_0K_0^{-2}/8)$
to be determined later. Thus,
$J_*=\log_2[R_0K_0^{-2}/(8C_*h)]\leq \log_2(2+1/h)$
and $J_*>1$ when $h<R_0K_0^{-2}/(16C_*) $.
Let
\begin{align*}
&Q_*(x_0)=\{(x,t)\in\Omega_T: \max (|x-x_0|,t^{1/2})\leq d_{J_*}\},
\\
&\Omega_*(x_0)=\{x\in \Omega: |x-x_0|\leq d_{J_*}\} \, ,  \\
&Q_j(x_0)=\{(x,t)\in \Omega_T:d_j\leq
\max (|x-x_0|,t^{1/2})\leq2d_j\},
\\
&\Omega_j(x_0)=\{x\in \Omega: d_j\leq|x-x_0|\leq2d_j\} ,
\\
&D_j(x_0)=\{x\in \Omega:  |x-x_0|\leq2d_j\}
\end{align*}
for $j\geq 1$; see Figure \ref{DomainFig}.

For $j=0$ we define $Q_0(x_0)=Q_T\backslash Q_1(x_0)$ and
$\Omega_0(x_0)=\Omega\backslash\Omega_1(x_0)$, and
for $j<0$ we simplify define $Q_j(x_0)=\Omega_j(x_0)=\emptyset$.
For all $j\geq 1$ we define
\begin{align*}
&\Omega_j'(x_0)=\Omega_{j-1}(x_0)\cup\Omega_{j}(x_0)\cup\Omega_{j+1}(x_0), \\
&\Omega_j''(x_0)=\Omega_{j-2}(x_0)\cup\Omega_{j}'(x_0)\cup\Omega_{j+2}(x_0),\\
&\Omega_j'''(x_0)=\Omega_{j-2}(x_0)\cup\Omega_{j}''(x_0)\cup\Omega_{j+2}(x_0),\\
& Q_j'(x_0)=Q_{j-1}(x_0)\cup Q_{j}(x_0)\cup Q_{j+1}(x_0), \\
& 
Q_j''(x_0)=Q_{j-2}(x_0)\cup Q_{j}'(x_0)\cup Q_{j+2}(x_0),\\
& 
Q_j'''(x_0)=Q_{j-2}(x_0)\cup Q_{j}''(x_0)\cup Q_{j+2}(x_0) ,\\
&D_j'(x_0)=D_{j-1}(x_0)\cup D_{j}(x_0) ,\\
&  D_j''(x_0)=D_{j-2}(x_0)\cup D_{j}'(x_0) ,\\
&  D_j'''(x_0)=D_{j-3}(x_0)\cup D_{j}''(x_0) .
\end{align*}
Then we have
\begin{align*}
&\Omega_T=\bigcup^{J_*}_{j=0}Q_j(x_0)\,\cup Q_*(x_0)
\quad\mbox{and}\quad
\Omega=\bigcup^{J_*}_{j=0}\Omega_j(x_0)\,\cup \Omega_*(x_0),
\end{align*}
We refer to $Q_*(x_0)$ as the ``innermost" set.
We shall write $\sum_{*,j}$ when the innermost set is included and
$\sum_j$ when it is not. When $x_0$ is fixed, if there is no ambiguity,
we simply write
$Q_j=Q_j(x_0)$, $Q_j'=Q_j'(x_0)$, $Q_j''=Q_j''(x_0)$,
$\Omega_j=\Omega_j(x_0)$, $\Omega_j'=\Omega_j'(x_0)$ and
$\Omega_j''=\Omega_j''(x_0)$.

In the rest of this paper, we denote by $C$ a generic positive constant, which
will be independent of $h$, $x_0$, and the undetermined constant
$C_*$ until it is determined at the end of section \ref{dka7}.

\subsection{Proof of Theorem \ref{MainTHM1}}
The keys to the proof of Theorem \ref{MainTHM1} are
several more precise estimates of the Green's functions.
Let $F(t)=\Gamma_h(t)-\Gamma(t)$.
Then for any $x_0 \in \Omega$, we have
\begin{align}
(E_h(t)v_h)(x_0) &= (F (t), v_h) + (\Gamma(t), v_h)
\nn\\
&=\int_0^t(\partial_tF (s), v_h)\d s+(F(0), v_h) + (\Gamma(t), v_h) ,
\label{EhtFt1} \\[10pt]
(t\partial_tE_h(t)v_h)(x_0)
&= (t\partial_tF (t), v_h) + (t\partial_t\Gamma(t), v_h)
\nn \\ &
=\int_0^t(s\partial_{ss}F (s)+\partial_{s}F (s), v_h)\d s
+ (t\partial_t\Gamma(t), v_h) ,
\label{EhtFt2}
\end{align}
with $\|F(0)\|_{L^1}=\|\widetilde\delta_{x_0}
-P_h\widetilde\delta_{x_0}\|_{L^1}\leq C $ (according to (P3) and (P4)).
Moreover, by the analyticity of the continuous parabolic semigroup
on $L^1(\Omega)$, we have
$$
\|\Gamma(t)\|_{L^1}+t\|\partial_t\Gamma(t)\|_{L^1}
\leq C \|\Gamma(0)\|_{L^1}
=C \|\widetilde\delta_{x_0}\|_{L^1}
\leq C .
$$

We present some estimates of these Green's
functions in the following lemma.
The proof of the lemma is the major work of this paper and
will be given in the next two sections.

\begin{lemma}
\label{lemma2-1}
{\it Under the assumptions of Theorem \ref{MainTHM1},
we have
\begin{align}
&\int_0^\infty\int_\Omega\big( |
\partial_tF(t,x,x_0) \big|
+\big|t\partial_{tt}F(t,x,x_0 ) |\big)\d x\d t \leq C , \label{FFEst2}\\
&|\nabla\partial_tG(t,x,x_0)|\leq C\max(t^{1/2},|x-x_0|)^{-3-N}
\quad\mbox{for}\,\,\,\, (x,t)\in \Omega\times(0,1) .
\label{DxtGrEst}
\end{align}
}
\end{lemma}

The estimates in Lemma \ref{lemma2-1} were proved in \cite{STW98}
for parabolic equations with the Neumann boundary condition
and in \cite{TW2} for the Dirichlet boundary condition.
However, their proofs are only valid for smooth coefficients
and smooth domains (as clearly mentioned in
their papers).
Later, these estimates were proved in \cite{Li}
for parabolic equations in smooth domains
of arbitrary dimensions under the Neumann boundary condition with
Lipschitz continuous coefficients.
Here we are concerned with the problem
in a convex polyhedron in two or three dimensional spaces under the Dirichlet
boundary condition with
$a_{ij} \in W^{1,N+\alpha}$.
\vskip0.1in

\noindent{\it Proof of Theorem \ref{MainTHM1}: }
Firstly, from \refe{EhtFt1}-\refe{EhtFt2} we see that
(\ref{STLEst}) is a consequence of (\ref{FFEst2}).

Secondly, we can view $E_h(t)$ as an analytic semigroup on 
$L^q(\Omega)$, defined by 
$$
(E_h(t)v)(x_0):=\int_\Omega \Gamma_h(t,x,x_0)v(x)\d x ,
$$
whose generator is $A_hP_h$. 
From \cite[Theorem 4.2]{Weis1} and
\cite[Lemma 4.c]{Weis2} (with a duality argument
for the case $q\geq 2$) we know that the maximal $L^p$
regularity \refe{LpqSt3} holds if
the following maximal ergodic estimate holds:
\begin{align}\label{maxergd}
\bigg\|\sup_{t>0}\frac{1}{t}\int_0^t|E_h(s)|\, v \,\d s\bigg\|_{L^q}
\leq C\|v\|_{L^q} ,\quad \forall\, v\in L^q(\Omega) ,
\end{align}
where 
$$
(|E_h(s)|v) (x_0):=\int_\Omega |\Gamma_h(t,x,x_0)| v(x)\d x .
$$
Let $G_{\rm tr}(t,x,x_0)$ be a truncated Green's function
which is symmetric with
respect to $x$ and $x_0$ and satisfies
$G_{\rm tr}(t,x,x_0)=G(t,x,x_0)$ when $(x,t)$ is outside
$Q_*(x_0)$ (see \cite[Section 4.2]{Li} on its construction).
Then we have
(assuming that $\tau_0^h$ is the triangle/tetrahedron which
contains $x_0$)
\begin{align*}
&\iint_{[\Omega\times(0,\infty)]\backslash
Q_*(x_0)}|\partial_t\Gamma(t,x,x_0)
-\partial_tG_{\rm tr}(t,x,x_0)|\d x\d t \\
&= \iint_{[\Omega\times(0,1)]\backslash Q_*(x_0)}\biggl|\int_\Omega
\partial_tG(t,x,y)\widetilde\delta_{x_0}(y)\d
y-\partial_tG(t,x,x_0)\biggl|\d x\d t\\
&\quad
+\iint_{\Omega\times(1,\infty)}|\partial_t\Gamma(t,x,x_0)
-\partial_tG(t,x,x_0)|\d x\d t\\
&\leq
Ch\iint_{[\Omega\times(0,1)]\backslash Q_*(x_0)}\sup_{y\in\tau_0^h}\big|
\nabla_y\partial_tG(t,x,y)\big|\d x\d t \\
&\quad
+C\int_{1}^\infty t^{-1}( \|\Gamma(t/2,\cdot,x_0)\|_{L^1}
+\|G(t/2,\cdot,x_0)\|_{L^1})\d t
\quad \mbox{[by semigroup estimate]} \\
&=Ch\sum_{j}\iint_{Q_j(x_0)}\sup_{(y,t)\in Q_j'(x)}\big|
\nabla_y\partial_tG(t,x,y)\big|\d x\d t
+C\int_{1}^\infty t^{-1-N/2}\d t
\quad \mbox{[see \refe{FEstP}] }\\
&\leq C\sum_{j}\frac{h}{d_j} +C
\quad \mbox{[see \refe{DxtGrEst}]}\\
&\leq C .
\end{align*}
By using energy estimates, it is easy to see
\begin{align*}
&\iint_{Q_*(x_0)}(|\partial_t\Gamma(t,x,x_0)|
+|\partial_tG_{\rm tr}(t,x,x_0)|)\d x\d t \\
&\leq
d_{J_*}^{N/2+1}
(\|\partial_t\Gamma(\cdot,\cdot,x_0)\|_{L^2(\Omega\times(0,1))}
+\|\partial_tG_{\rm tr}(\cdot,\cdot,x_0)\|_{L^2(\Omega\times(0,1))})\leq
C_*^{N/2+1} ,
\end{align*}
where the constant $C_*$ will be determined
at the end of Section 4.
Then \refe{FFEst2} and the last two inequalities imply
\begin{align*}
&\int_0^\infty\int_{\Omega}|\partial_t\Gamma_h(t,x,x_0)
-\partial_tG_{\rm tr}(t,x,x_0)|\d x\d t \leq C .
\end{align*}
In other words, the symmetric kernel
$K(x,y):=\int_0^\infty|\partial_t\Gamma_h(t,x,y)
-\partial_tG^*_{\rm tr}(t,x,y)|\d t$
satisfies
$$
\sup_{y\in\Omega}\int_\Omega K(x,y)\d
x+\sup_{x\in\Omega}\int_\Omega K(x,y)\d y\leq C ,
$$
and therefore, Schur's lemma implies that the corresponding operator $M_K$,
defined by $M_K{v}(x)=\int_\Omega K(x,y){v}(y)\d y$, is bounded on
$L^q(\Omega)$ for all $1\leq q\leq\infty$.
Let
$E^*_{\rm tr}(t){v}(x)=\int_\Omega G_{\rm tr}^*(t,x,y){v}(y)\d y$
and note that
$E^*_{\rm tr}(t){v}(x)\leq E(t)|{v}|(x)$ (because
$G_{\rm tr}^*(t,x,y)\leq G(t,x,y)$).
We have
\begin{align*}
&{ \sup_{t>0}(|E_h(t)P_h|v)(x)}\\
& 
\leq 
{ \sup_{t>0}(|E_h(t)P_h| v - E^*_{\rm tr}(t) v)(x) 
+\sup_{t>0}(E^*_{\rm tr}(t) v)(x) }\\
& 
\leq { \sup_{t>0}\big(|E_h(t)P_h-E^*_{\rm
tr}(t)|\,|v|\big)(x)+\sup_{t>0}(E^*_{\rm tr}(t)v)(x) }\\
& 
={
\sup_{t>0}\bigg|(|P_h\widetilde\delta_{x}|,|v|)+\int_0^t
\int_\Omega|
\partial_t\Gamma_h(s,x,y)
-\partial_tG^*_{\rm tr}(s,x,y)|\, |v(y)|\d
y\d s\bigg| }\\
&\quad +
{ \sup_{t>0}(E^*_{\rm tr}(t)|v|)(x) }\\
&
\leq 
{ (|P_h\widetilde\delta_{x}|,|v|)
+(M_K|v|)(x)+\sup_{t>0}(E(t)|v|)(x)
}
\end{align*} 
where
$$
{ \|\sup_{t>0}E(t) |v|\|_{L^q}\leq C_q\|v\|_{L^q}, }
\quad\forall~1<q<\infty ,
$$
is a simple consequence of the Gaussian estimate \refe{FEstP}
(Corollary 2.1.12 and Theorem 2.1.6 of \cite{Grafakos}).
This proves a stronger estimate than
\refe{maxergd}.  The proof of \refe{LpqSt3} is completed.

Finally, \refe{PDE0}-\refe{FEEq0}
imply that the error $e_h=P_hu-u_h$ satisfies the equation
(when $u^0=u^0_h=0$)
\begin{align}
\partial_t(A_h^{-1}e_h)+A_h(A_h^{-1}e_h)=P_hu-R_hu .
\end{align}
By applying \refe{LpqSt3} to the equation above, we obtain
\begin{align}
\|e_h\|_{L^p((0,T);L^q)}
\leq C_{p,q}\|P_hu-R_hu\|_{L^p((0,T);L^q)}
\leq C_{p,q}h\|u\|_{L^p((0,T);W^{1,q})}
\end{align}
for $1<p<\infty$ and $2\leq q<\infty$,
where we have used the inequality
$\|P_hu-R_hu\|_{L^q}\leq  C_{q}h\|u\|_{W^{1,q}}$,
which only holds for $2\leq q<\infty$
in convex polygons/polyhedra.
Then, by using an inverse inequality
and \refe{LpW1quf}, we have
\begin{align*}
\|e_h\|_{L^p((0,T);W^{1,q})}
&\leq Ch^{-1}\|e_h\|_{L^p((0,T);L^q)} \\
&\leq C_{p,q}\|u\|_{L^p((0,T);W^{1,q})} \\
&\leq C_{p,q}\|f\|_{L^p((0,T);W^{-1,q})} ,
\end{align*}
which implies \refe{LpqSt2} for the case
$1<p<\infty$ and $2\leq q<\infty$.

In the case $1<p<\infty$ and $1<q\leq 2$,
we define $\vec g=\nabla \Delta^{-1}P_hf$ and express
the solution of \refe{FEEq0} by
(when $u_h^0=0$)
\begin{align*}
\nabla u_h=
{\mathcal L}_h \vec g:=\int_0^t\nabla A_h^{-1/2} 
A_hE_h(t-s) A_h^{-1/2}\nabla\cdot
\vec g(s) \d s .
\end{align*}
In order to prove the boundedness of the operator
${\mathcal L}_h$ on $L^p((0,T);(L^q)^N)$, we only need to
prove the boundedness of its dual operator
${\mathcal L}_h'$ on $L^{p'}((0,T);(L^{q'})^N)$.
It is easy to see that
\begin{align*}
\int_0^T({\mathcal L}_h \vec g , \vec \eta)\d t
=\int_0^T\bigg(\vec g,\int_s^T
\nabla A_h^{-1/2} A_hE_h(t-s) A_h^{-1/2}\nabla\cdot
\vec \eta(t) \d t\bigg)\d s ,
\end{align*}
which gives
\begin{align*}
{\mathcal L}_h' \vec \eta:=\int_s^T\nabla A_h^{-1/2}
A_hE_h(t-s) A_h^{-1/2}\nabla\cdot
\vec \eta(s) \d s .
\end{align*}
If we define the backward finite element problem
\begin{align}
\left\{
\begin{array}{ll}
-\big(\partial_tw_h,v_h\big)
+ (a\nabla w_h, \nabla v_h)
=(\nabla\cdot\vec\eta,v_h) , \,\,\,\,\forall\, v_h\in S_h,
\\
w_h(T)=0 ,
\end{array}
\right.
\end{align}
then ${\mathcal L}_h' \vec \eta=\nabla w_h$.
By a time reversal we obtain, as shown in the last paragraph,
\begin{align*}
\|\nabla w_h\|_{L^{p'}((0,T);L^{q'})}
\leq C_{p,q}\|\nabla\cdot\vec\eta\|_{L^{p'}((0,T);W^{-1,q'})}
\leq C_{p,q}\|\vec \eta\|_{L^{p'}((0,T);L^{q'})}  ,
\end{align*}
for $1<p'<\infty$ and $2\leq q'<\infty$,
which implies the boundedness of
${\mathcal L}_h'$ on $L^{p'}((0,T);(L^{q'})^N)$.
By duality, we derive the boundedness of
${\mathcal L}_h$ on $L^p((0,T);(L^q)^N)$ and therefore,
\begin{align*}
\|\nabla u_h\|_{L^{p}((0,T);L^{q})}
&\leq C_{p,q}\|\vec g\|_{L^{p}((0,T);L^{q})} \\
&\leq C_{p,q}\|P_hf\|_{L^p(0,T);W^{-1,q})} \\
&\leq C_{p,q}\|f\|_{L^p(0,T);W^{-1,q})}   .
\end{align*}
This proves \refe{LpqSt2} in the case $1<p<\infty$ and $1<q\leq 2$.

The proof of Theorem \ref{MainTHM1} is completed
(based on Lemma \ref{lemma2-1}).
\,\endproof \medskip

\noindent{\bf Remark 2.1}$\,\,$
In the proof of \refe{LpqSt2}, we have used
an $L^q$ error estimate of the Ritz projection for $2\leq q<\infty$, which
can be proved in the same way as used in \cite[Corollary]{RS} by using
the $W^{1,q}$-stability of the Ritz projection.
This $W^{1,q}$-stability
is based on an interpolation between these two cases
$q=2$ and $q=\infty$.
The case $q=2$ is trivial and the
case $q=\infty$ was studied by several authors,
such as \cite{RS} for $r=1$ and 2D convex polygons 
(which requires $H^2$ regularity of the elliptic problem),
\cite{Schatz} for $r\geq 2$ and 2D arbitrary polygons 
(as a consequence of the $L^\infty$ stability proved therein,
which only requires $H^{3/2+\varepsilon}$ regularity
of the elliptic problem),
and \cite{GLRS} for $r\geq 1$ in 3D convex polyhedra
(which requires $H^2$ and $C^{1+\alpha}$ regularity of the elliptic problem).
These essential properties used by \cite{GLRS,RS,Schatz} are all
possessed by the elliptic problem when
the domain is convex polygonal/polyhedral and the coefficients
$a_{ij}$ are $W^{1,N+\alpha}$.

In the rest of this paper, we focus on the proof of Lemma
\ref{lemma2-1}.

\section{Preliminary analysis}
\setcounter{equation}{0}
In this section,  we present two propositions.

\subsection{Superapproximation of smoothly
truncated finite element functions}
In this subsection, we prove the following proposition,
which is needed in proving Lemma \ref{lemma2-1}.
\begin{proposition}\label{lemsupa0}
{\it If $0\leq \omega \leq 1$ is a smooth cut-off function which equals zero in
$\Omega\backslash D$, satisfying
$|\partial^\beta\omega |\leq Cd^{-|\beta|}$
for all multi-index $\beta$ such that
$|\beta|=0,1,\cdots,r+1$ and $d\geq 10\kappa h$,
then for any $\psi_h\in S_h$ there exists $\chi_h\in S_h^0(B_d(D))$
such that
\begin{align}\label{lemsupa}
d^2\|R_h(\omega\psi_h)-\chi_h\|_{H^1}
+d\|\omega \psi_h-\chi_h\|_{L^2}\leq Ch\|\psi_h\|_{L^2(B_d(D))} .
\end{align}
}
\end{proposition}
\noindent{\it Proof.}$\,\,$
Define $0\leq \widetilde\omega\leq 1$ as
a smooth cut-off function which is zero
outside $B_{0.8d}(D)$, satisfying that $\widetilde \omega=1$
on $B_{0.7d}(D)$ and
$|\partial^\beta\widetilde\omega |\leq Cd^{-|\beta|}$
for $|\beta|=0,1,\cdots,r+1$.

First we prove the following inequality
\begin{align}\label{PInEq1}
\|\omega\psi_h-R_h(\omega\psi_h) \|_{L^2 }
\leq Chd^{-1}\|\psi_h\|_{L^2(B_{0.3d}(D))}
\end{align}
by a duality argument. We define
$\phi$ as the solution of the elliptic PDE
$$
\left\{
\begin{array}{ll}
-\nabla\cdot(a\nabla \phi)= v &\mbox{in}\,\,\,\, \Omega,\\
\phi=0 &\mbox{on}\,\,\, \partial\Omega  .
\end{array}\right.
$$
We see that
\begin{align*}
\big( v  ,\omega\psi_h-R_h(\omega\psi_h)\big) &=
\big(a\nabla \phi,\nabla (\omega\psi_h-R_h(\omega\psi_h))\big)  \\
&=\big(a\nabla (\phi -R_h\phi),\nabla (\omega\psi_h-R_h(\omega\psi_h))\big)\\
&=\big(a\nabla (\phi -R_h\phi),\nabla (\omega\psi_h-\Pi_h(\omega\psi_h))\big ) \\
&\leq C\|\phi -R_h\phi \|_{H^1}\|\omega\psi_h-\Pi_h(\omega\psi_h)\|_{H^1} \\
&\leq Chd^{-1}\|v\|_{L^2}  \|\psi_h\|_{L^2(B_{0.3d}(D))} .
\end{align*}
where we have used the superapproximation property
(P2) in section \ref{Sec2-2}
and the $H^1$ error estimate:
$$
\|\phi -R_h\phi \|_{H^1} \leq C h\|\phi\|_{H^2}
\leq Ch\|v\|_{L^2} .
$$
\refe{PInEq1} follows these inequalities.

Secondly, it is noted that the following inequality
\begin{align}\label{PInEq2}
\|R_h(\omega\psi_h) \|_{H^1(B_{d}(D)\backslash B_{0.5d}(D))}
\leq Cd^{-1}\|R_h(\omega\psi_h) \|_{L^2(B_{d}(D)\backslash B_{0.3d}(D))}
\end{align}
was proved in Lemma 4.4 of \cite{SW77}
(also see Page 1374 of \cite{STW98})
as a consequence of the discrete elliptic equation
$$
(a\nabla R_h(\omega\psi_h),\nabla\eta)=0 ,\quad\mbox{ for }\,\,\,\,
 \eta\in S_h^0(B_{d}(D)\backslash D)  \, .
$$

Let $\chi_h=\Pi_h[\widetilde\omega R_h(\omega\psi_h)]$
and note that the support of $\chi_h$
is contained in $B_{0.8d}(D)$.
By using the superapproximation
property (P2), we have
\begin{align}\label{PInEq3}
&d^{-1}\|R_h(\omega\psi_h)-\chi_h\|_{L^2}
+\|R_h(\omega\psi_h)-\chi_h\|_{H^1} \nn\\
&\leq d^{-1}\|R_h(\omega\psi_h)-\widetilde \omega R_h(\omega\psi_h) \|_{L^2}
+d^{-1}\|\widetilde \omega R_h(\omega\psi_h)
-\Pi_h[\widetilde \omega R_h(\omega\psi_h)]\|_{L^2} \nn \\
&\quad
+\|R_h(\omega\psi_h)-\widetilde \omega R_h(\omega\psi_h) \|_{H^1}
+\|\widetilde \omega R_h(\omega\psi_h)
-\Pi_h[\widetilde \omega R_h(\omega\psi_h)]\|_{H^1} \nn \\
&\leq Cd^{-1}\|R_h(\omega\psi_h)\|_{L^2(B_{d}(D)\backslash B_{0.3d}(D))} \nn\\
&=Cd^{-1}\|R_h(\omega\psi_h)-\omega\psi_h \|_{L^2(B_{d}(D)\backslash B_{0.3d}(D))}
\,\,\quad\mbox{(because $\omega=0$ on
$B_{d}(D)\backslash B_{0.3d}(D)$)} \nn\\
&\leq Chd^{-2}\|\psi_h\|_{L^2(B_{0.3d}(D))}
\qquad\qquad\qquad\qquad\quad
\,\mbox{(as a consequence of \refe{PInEq1})}
\end{align}
and from \refe{PInEq1} we see that
\begin{align}
\|\omega\psi_h -\chi_h\|_{L^2}
&\leq \|\omega\psi_h -R_h(\omega\psi_h)\|_{L^2}
+\|R_h(\omega\psi_h) -\chi_h\|_{L^2} \nn\\
&\leq Chd^{-1}\|\psi_h\|_{L^2(B_{0.3d}(D))} \, .
\end{align}
\refe{lemsupa} follows immediately and the proof of
the Proposition \ref{lemsupa0} is completed.
\quad \endproof \medskip

\noindent{\bf Remark 3.1}$\,\,\,$
In the proof of Proposition \ref{lemsupa0}
we have assumed that $d\geq 10\kappa d$ and used
$B_{0.3d}(D)$, $B_{0.7d}(D)$, $B_{0.8d}(D)$ ...
to make sure that their radius differ from each other
by at least $\kappa h$ so that
the superapproximation property (P2) can be used.

\subsection{Local error estimate}
The following proposition is
concerned with a local energy
error estimate of parabolic equations.

\begin{proposition}\label{LocEEst}
{\it
Suppose that $\phi,\partial_t\phi\in L^2((0,T);H^1_0)$ and
$\phi_h\in H^1((0,T);S_h)$, and
$e=\phi_h-\phi$ satisfies the equation
\begin{align}
(e_t,\chi)+(a\nabla e,\nabla \chi) =0,
\quad\forall\, \chi\in S_h,~t>0 ,
\label{p32}
\end{align}
with $\phi(0)=0$ in $\Omega_j''''$.
Then for any $m>0$,
there exists a constant $C_m>0$,
independent of $h$ and $d_j$, such that
\begin{align}
&\|e_t\|_{L^2(Q_j)} + d_j^{-1}\|\nabla e\|_{L^2(Q_j)} \nn\\
&\leq
C_m\big(I_j(\phi_{h}(0))+X_j(\Pi_h\phi-\phi)
+H_j(e)+d_j^{-2}\|e\|_{L^2(Q_j''')}\big),
\label{3-1}
\end{align}
where
\begin{align*}
&I_j(\phi_{h}(0))=d_j^{-1}\|\phi_{h}(0)\|_{L^2(\Omega_j''')}
+\|\phi_{h}(0)\|_{H^1(\Omega_j''')} ,\\
&X_j(\Pi_h\phi-\phi)=d_j\|\nabla\partial_t(\Pi_h\phi-\phi)\|_{L^2(Q_j''')}
+\|\partial_t(\Pi_h\phi-\phi)\|_{L^2(Q_j''')} \\
&\qquad \qquad \qquad\
+d_j^{-1}\|\nabla(\Pi_h\phi-\phi)\|_{L^2(Q_j''')}+
d_j^{-2}\|\Pi_h\phi-\phi\|_{L^2(Q_j''')},\\
&H_j(e)=(h/d_j)^m\big(\|e_t\|_{L^2(Q_j''')}
+d_j^{-1}\|\nabla e\|_{L^2(Q_j''')}\big) .
\end{align*}
}
\end{proposition}

Before we prove Proposition \ref{LocEEst},
we present a local energy estimate for
finite element solutions of parabolic equations.
\begin{lemma}\label{lemlocEng0}
{\it
Suppose that $\phi_h(t)\in S_h$ satisfies
\begin{align*}
&\big(\partial_t\phi_h ,\chi\big)
+\big(a\nabla \phi_h ,\nabla\chi\big) =0,
\qquad\mbox{for}\,\,\,\,\chi\in S_h^0(\Omega_j''),\,\,\, t\in(0,d_j^2]  ,\\
&\big(\partial_t \phi_h ,\chi\big)
+\big(a\nabla \phi_h ,\nabla\chi\big)
=0   ,\qquad\mbox{for}\,\,\,\,\chi\in S_h^0(D_j'') ,\,\,\, t\in(d_j^2/4,4d_j^2) .
\end{align*}
Then for any $m>0$ there exists $C_m>0$,
independent of $h$ and $d_j$, such that
\begin{align}\label{locEng0}
&\|\partial_t\phi_h\|_{L^2(Q_j)} +d_j^{-1}\|\nabla \phi_h\|_{L^2(Q_j)} \nn\\
&\leq  C_m\left(\|\nabla \phi_h(0)\|_{L^2(\Omega_j'')}
+d_j^{-1}\|\phi_h(0)\|_{L^2(\Omega_j'')} \right) \nn\\
&\quad
+ C_m\left(\frac{h}{d_j}\right)^{m}\big(
\|\partial_t\phi_h\|_{L^2(Q_{j}'')}
+ d_j^{-1}\|\nabla\phi_h\|_{L^2(Q_{j}'')}\big)
+ C_md_j^{-2} \|\phi_h\|_{L^2(Q_{j}'')} .
\end{align}
}
\end{lemma}

\noindent{\it Proof.}$\,\,$
Note that $Q_j=[\Omega_j\times(0,d_j^2) ]
\cup [D_j\times (d_j^2,4d_j^2)]$.
We first present estimates in the domain $\Omega_j\times(0,d_j^2)$
and then present estimates in the domain $D_j\times (d_j^2,4d_j^2)$.

Let $\omega$ be a smooth cut-off function
which equals $1$ in $\Omega_j$
and vanishes outside $\Omega_j'$,
and let $\widetilde\omega$ be a smooth cut-off function
which equals $1$ in $\Omega_j'''$
and vanishes outside $\Omega_j''''$ with
\begin{align}
|\partial^{\beta}\omega|\leq Cd_j^{-|\beta|}
\qquad\mbox{and}\qquad
|\partial^{\beta}\widetilde\omega|\leq Cd_j^{-|\beta|}.
\label{3-2}
\end{align}
Let $v_h=\Pi_h(\widetilde\omega\phi_h)\in S_h^0(\Omega_j''''')$
so that $v_h=\phi_h$ in $\Omega_j''$, satisfying
(due to the superapproximation property (P2))
\begin{align*}
&\|v_h\|_{L^2 }
\leq C\|\phi_h\|_{L^2(\Omega_j'''')}  ,\\
&\|\nabla v_h\|_{L^2 }
\leq C\|\nabla \phi_h\|_{L^2(\Omega_j'''')}
+ Cd_j^{-1}\| \phi_h\|_{L^2(\Omega_j'''')}
\end{align*}
and
$$
\big(\partial_tv_h,\chi\big)+\big(a\nabla v_h,\nabla\chi\big)
=0  ,\quad\forall\, \chi\in S_h^0(\Omega_j'') \, .
$$
It follows that
\begin{align*}
&\frac{1}{2}\frac{\d}{\d t}\|\omega v_h\|^2
+(\omega^2 a\nabla v_h,\nabla v_h)\\
&=\big[\big(\partial_tv_h,\omega^2v_h-\chi_h\big)
+\big(a\nabla v_h,\nabla(R_h(\omega^2v_h)-\chi_h)\big)\big]\\
&\quad 
+\big[\big(\omega^2a\nabla v_h,\nabla v_h\big)
-\big(a\nabla v_h,\nabla(\omega^2v_h)\big)\big] \\
&\leq
\big[ C\|\partial_tv_h\|_{L^2}\|v_h\|_{L^2}hd_j^{-1}
+ C\|\nabla v_h\|_{L^2}\|v_h\|_{L^2}hd_j^{-2}  \big]
+ C\big(  \omega a\nabla v_h,2 v_h\nabla \omega\big),
\end{align*}
where we have used (P2) and Proposition \ref{lemsupa0},
and from \refe{3-2} we see that
\begin{align*}
\big(  \omega a\nabla v_h,2 v_h\nabla \omega\big)
\leq \big(  |\omega a\nabla v_h|, \, 2 |v_h| \big) d_j^{-1}
\leq C \|\omega a\nabla v_h \|_{L^2} \|v_h\|_{L^2} d_j^{-1}
\, .
\end{align*}
The last two inequalities imply
\begin{align}\label{locEng1}
&\|\phi_h\|_{L^\infty((0,d_j^2)\times L^2(\Omega_j))}
+\|\nabla \phi_h\|_{L^2((0,d_j^2);L^2(\Omega_j))} \nn\\
&\leq C\|\phi_h(0)\|_{L^2(\Omega_j'''')}
+C\|\partial_t\phi_h\|_{L^2((0,d_j^2);L^2(\Omega_j''''))}h \nn \\
&\quad
+C\|\nabla \phi_h\|_{L^2((0,d_j^2);L^2(\Omega_j''''))} h d_j^{-1}
+C\|\phi_h\|_{L^2((0,d_j^2);L^2(\Omega_j''''))} d_j^{-1} .
\end{align}

By using Proposition \ref{lemsupa0} again, we derive that
\begin{align*}
&\|\omega^2\partial_tv_h\|_{L^2}^2
+\frac{1}{2}\frac{\d}{\d t}(\omega^4 a\nabla v_h,\nabla v_h) \\
&=\big[(\partial_tv_h,\omega^4\partial_tv_h-\chi_h)
+(a\nabla v_h,\nabla[R_h(\omega^4\partial_tv_h)-\chi_h])\big]
+(4\omega^3a\nabla v_h,\partial_tv_h\nabla \omega)  \\
&\leq C\|\partial_tv_h\|_{L^2}^2hd_j^{-1}
+C\|\nabla v_h\|_{L^2}\|\partial_tv_h\|_{L^2}hd_j^{-2}
+C\|\omega\nabla v_h\|_{L^2}
\|\omega^2\partial_t v_h\|_{L^2} d_j^{-1}  \\
&\leq C\|\partial_tv_h\|_{L^2}^2hd_j^{-1}
+ C\|\nabla v_h\|_{L^2}^2 d_j^{-2}
+\frac{1}{2} \|\omega^2\partial_t v_h\|_{L^2}^2 ,
\end{align*}
which reduces to
\begin{align*}
&\|\omega^2\partial_tv_h\|_{L^2((0,d_j^2);L^2(\Omega))}^2 \\
&\leq C\|\nabla v_h(0)\|_{L^2(\Omega)}^2
+ C\|\partial_tv_h\|_{L^2((0,d_j^2);L^2(\Omega))}^2hd_j^{-1}
+ C\|\nabla v_h\|_{L^2((0,d_j^2);L^2(\Omega))}^2 d_j^{-2}  .
\end{align*}
The inequality above further implies
\begin{align}\label{locEng2}
\|\partial_t\phi_h\|_{L^2((0,d_j^2);L^2(\Omega_j))}
&\leq C(\|\nabla\phi_h(0)\|_{L^2(\Omega_j'''')}
+d_j^{-1}\|\phi_h(0)\|_{L^2(\Omega_j'''')}) \nn\\
&\quad
+ C\|\partial_t\phi_h\|_{L^2((0,d_j^2);L^2(\Omega_j''''))} h^{1/2}d_j^{-1/2}
\nn\\
&\quad
+ Cd_j^{-1}(\|\phi_h\|_{L^\infty((0,d_j^2)\times L^2(\Omega_j''''))}
+\|\nabla \phi_h\|_{L^2((0,d_j^2);L^2(\Omega_j''''))} )  .
\end{align}
With an obvious change of indices
(from $\Omega''''$ to $\Omega'''$ on the right-hand side,
and from $\Omega$ to $\Omega'$ on the left-hand side),
\refe{locEng1}-\refe{locEng2}
imply
\begin{align}\label{locEng1-2}
& \|\phi_h\|_{L^\infty((0,d_j^2); L^2(\Omega_j'))}
+\|\nabla \phi_h\|_{L^2((0,d_j^2);L^2(\Omega_j'))} \nn \\
&\leq C\|\phi_h(0)\|_{L^2(\Omega_j''')}
+C h \|\partial_t\phi_h\|_{L^2((0,d_j^2);L^2(\Omega_j'''))} \nn \\
&\quad
+C h d_j^{-1} \|\nabla \phi_h\|_{L^2((0,d_j^2);L^2(\Omega_j'''))}
+C d_j^{-1} \|\phi_h\|_{L^2((0,d_j^2);L^2(\Omega_j'''))} .
\end{align}
and
\begin{align}\label{locEng2-2}
\|\partial_t\phi_h\|_{L^2((0,d_j^2);L^2(\Omega_j))}
&\leq C(\|\nabla\phi_h(0)\|_{L^2(\Omega_j')}
+d_j^{-1}\|\phi_h(0)\|_{L^2(\Omega_j')}) \nn\\
&\quad
+ C h^{1/2}d_j^{-1/2} \|\partial_t\phi_h\|_{L^2((0,d_j^2);L^2(\Omega_j'))}
\nn\\
&\quad
+ Cd_j^{-1} (\|\phi_h\|_{L^\infty((0,d_j^2); L^2(\Omega_j'))}
+\|\nabla \phi_h\|_{L^2((0,d_j^2);L^2(\Omega_j'))} )  .
\end{align}

In the same way as we derive \refe{locEng1-2}-\refe{locEng2-2},
by choosing $\overline\omega(x,t)=\omega_1(x)\omega_2(t)$
with $\omega_1=1$ in $D_j'$,  $\omega_1=0$ outside $D_j''$,
$\omega_2=1$ for $t\in(d_j^2,4d_j^2)$ and
$\omega_2=0$ for $t\in(0,d_j^2/2)$,
we can derive that
\begin{align}\label{locEng3}
&\|\overline\omega\phi_h\|_{L^\infty((0,4d_j^2); L^2(\Omega))}
+\|\overline\omega\nabla \phi_h\|_{L^2((0,4d_j^2);L^2(\Omega))} \nn\\
&\leq  C h \|\partial_t\phi_h\|_{L^2((d_j^2/4,4d_j^2);L^2(D_{j}'''))} \nn\\
&\quad
+C h d_j^{-1} \|\nabla \phi_h\|_{L^2((d_j^2/4,4d_j^2);L^2(D_{j}'''))} \nn\\
&\quad
+C d_j^{-1} \|\phi_h\|_{L^2((d_j^2/4,4d_j^2);L^2(D_{j}'''))}
\end{align}
and
\begin{align}\label{locEng4}
\|\partial_t \phi_h\|_{L^2((d_j^2,4d_j^2);L^2(D_{j}))}
&\le C h^{1/2}d_j^{-1/2}
\|\partial_t \phi_h\|_{L^2((d_j^2/4,4d_j^2);L^2(D_{j}'''))}
\nn\\
&\quad
+ Cd_j^{-1} (\|\overline\omega\phi_h\|_{L^\infty((0,4d_j^2); L^2(\Omega))}
+\|\overline\omega\nabla  \phi_h\|_{L^2((0,4d_j^2);L^2(\Omega))}  )  .
\end{align}
By noting the definition of $\omega$ and $\overline \omega$,
we have
\begin{align*}
& \|\partial_t \phi_h\|_{L^2(Q_j)} \le
C ( \|\partial_t\phi_h\|_{L^2((0,d_j^2);L^2(\Omega_j))}
+ \|\partial_t \phi_h\|_{L^2((d_j^2,4d_j^2);L^2(D_{j}))} )
\\
& \|\nabla  \phi_h\|_{L^2(Q_j)} \le
C( \|\omega\nabla \phi_h\|_{L^2((0,d_j^2);L^2(\Omega))}
+ \|\overline\omega\nabla \phi_h\|_{L^2((0,4d_j^2);L^2(\Omega))} )
\, .
\end{align*}
With the last two inequalities, combining \refe{locEng1}-\refe{locEng4} gives
\begin{align*}
\|\partial_t \phi_h\|_{L^2(Q_j)} +d_j^{-1}\|\nabla  \phi_h\|_{L^2(Q_j)}
&\leq  C(\|\nabla  \phi_h(0)\|_{L^2(\Omega_j''')}
+d_j^{-1}\| \phi_h(0)\|_{L^2(\Omega_j''')} )  \\
&\quad
+ C\left(\frac{h}{d_j}\right)^{\frac{1}{2}}
\big( \|\partial_t \phi_h\|_{L^2(Q_{j}''' )}
+ d_j^{-1}\|\nabla\phi_h\|_{L^2(Q_{j}''' )}\big) \\
&\quad
+ Cd_j^{-2} \|  \phi_h\|_{L^2(Q_{j}''')}    .
\end{align*}
Iterating the inequality above and
changing the indices, we derive \refe{locEng0}.
\quad \endproof
\vskip0.1in

Now we are ready to prove Proposition \ref{LocEEst}.
\medskip

\noindent{\it Proof of Proposition \ref{LocEEst}}$\quad$
Let $\widetilde\omega(x,t)$ be a smooth cut-off function
which equals $1$ in $Q_j''$
and vanishes outside $Q_j'''$,
and let $\widetilde  \phi=\widetilde\omega  \phi$.
Then $\widetilde\phi(0)=0$ and we have
\begin{align*}
&\big(\partial_t(\widetilde \phi-\phi_h),\chi\big)
+\big(a\nabla (\widetilde \phi-\phi_h),\nabla\chi\big)
=0  ,\qquad\mbox{for}\,\,\,\,\chi\in S_h^0(\Omega_j ') ,\,\,\, t\in(0,d_j^2] ,\\
& \big(\partial_t(\widetilde \phi-\phi_h),\chi\big)
+\big(a\nabla (\widetilde \phi-\phi_h),\nabla\chi\big)
=0  ,\qquad\mbox{for}\,\,\,\,\chi\in S_h^0(D_j ') ,\,\,\, t\in(d_j^2/4,4d_j^2) .
\end{align*}
Let $\widetilde \phi_h\in S_h$ be the solution of
\begin{align}\label{wdtzh}
\big(\partial_t(\widetilde\phi-\widetilde \phi_h),\chi_h\big)
+\big(a\nabla (\widetilde\phi-\widetilde \phi_h),\nabla\chi_h\big)
=0   ,\qquad\mbox{for}\,\,\,\,\chi_h\in S_h
\end{align}
with $\widetilde\phi_h(0)=\Pi_h\widetilde \phi(0)= 0$ so that
\begin{align}
&\big(\partial_t(\widetilde \phi_h-\phi_h),\chi_h\big)
+\big(a\nabla (\widetilde \phi_h-\phi_h),\nabla\chi_h\big)
=0   ,\quad\forall\,\,\chi_h\in S_h^0(\Omega_j'),\,\,\, t\in(0,d_j^2]  ,
\label{wdtzh2}\\
&\big(\partial_t(\widetilde \phi_h-\phi_h),\chi_h\big)
+\big(a\nabla (\widetilde \phi_h- \phi_h),\nabla\chi_h\big)
=0   ,\quad\forall\,\,\chi_h\in S_h^0(D_j') ,\,\,\, t\in(d_j^2/4,4d_j^2) .
\label{wdtzh3}
\end{align}
Substituting $\chi_h=P_h\widetilde\phi-\widetilde\phi_h$
into (\ref{wdtzh}) we obtain
\begin{align*}
&\|\widetilde\phi-\widetilde \phi_h\|_{L^\infty((0,T);L^2)}^2
+\int_0^T\big(a\nabla (\widetilde\phi-\widetilde\phi_h),
\nabla(\widetilde\phi-\widetilde\phi_h)\big)\d t   \\
&=\int_0^T\big(\partial_t(\widetilde\phi-\widetilde \phi_h),
\widetilde\phi-P_h\widetilde\phi \big)\d t
+\int_0^T\big(a\nabla (\widetilde\phi-\widetilde\phi_h),
\nabla(\widetilde\phi-P_h\widetilde\phi)\big)\d t   ,
\end{align*}
which implies
\begin{align} \label{H1phiphih}
&\|\widetilde\phi-\widetilde \phi_h\|_{L^\infty((0,T); L^2(Q_T)}^2
+\|\nabla (\widetilde \phi-\widetilde \phi_h)\|_{L^2(Q_T)}^2 \nn\\
&
\leq C\|\partial_t (\widetilde \phi-\widetilde \phi_h)\|_{L^2(Q_T)}
\| \widetilde \phi \|_{L^2(Q_T)}
+C\| \nabla\widetilde\phi \|_{L^2(Q_T)}^2 .
\end{align}
Substituting $\chi_h=\partial_t(P_h\widetilde\phi-\widetilde\phi_h)$
into (\ref{wdtzh}) we obtain
\begin{align*}
&\|\partial_t(\widetilde\phi-\widetilde \phi_h)\|_{L^2(Q_T)}^2
+\sup_{0\leq t\leq T}\big(a\nabla (\widetilde\phi-\widetilde\phi_h),
\nabla(\widetilde\phi-\widetilde\phi_h)\big)   \\
&=\int_0^T\big(\partial_t(\widetilde\phi-\widetilde \phi_h),
\partial_t(\widetilde\phi-P_h\widetilde\phi)\big)\d t
+\int_0^T\big(a\nabla (\widetilde\phi-\widetilde\phi_h),
\nabla\partial_t(\widetilde\phi-P_h\widetilde\phi)\big)\d t  \\
&\leq \|\partial_t(\widetilde\phi-\widetilde \phi_h)\|_{L^2(Q_T)}
\|\partial_t\widetilde\phi \|_{L^2(Q_T)}
+\|\nabla(\widetilde\phi-\widetilde \phi_h)\|_{L^2(Q_T)}
\|\nabla\partial_t\widetilde\phi\|_{L^2(Q_T)}   ,
\end{align*}
which implies
\begin{align} \label{H1tphiphih}
&\|\partial_t(\widetilde \phi-\widetilde \phi_h)\|_{L^2(Q_T)}^2
\leq C\|\partial_t  \widetilde \phi \|_{L^2(Q_T)}^2
+C\|\nabla (\widetilde \phi-\widetilde \phi_h)\|_{L^2(Q_T)}
\| \nabla\partial_t\widetilde \phi \|_{L^2(Q_T)} ,
\end{align}
It follows that
\begin{align*}
&\|\partial_t(\widetilde \phi-\widetilde\phi_h)\|_{L^2(Q_T)}^2
+d_j^{-2}\|\nabla (\widetilde \phi-\widetilde\phi_h)\|_{L^2(Q_T)}^2
+d_j^{-2}\|\widetilde\phi-\widetilde \phi_h\|_{L^\infty((0,T);L^2)}^2\\
&\leq \frac{1}{2}\|\partial_t (\widetilde\phi-\widetilde\phi_h)\|_{L^2(Q_T)}^2
+Cd_j^{-4}\| \widetilde \phi \|_{L^2(Q_T)}^2
+Cd_j^{-2}\| \nabla\widetilde \phi \|_{L^2(Q_T)}^2,\\
&\quad
+C\|\partial_t  \widetilde \phi \|_{L^2(Q_T)}^2
+\frac{1}{2}d_j^{-2}\|\nabla (\widetilde \phi-\widetilde \phi_h)\|_{L^2(Q_T)}^2
+Cd_j^{2}\| \nabla\partial_t\widetilde\phi \|_{L^2(Q_T)}^2 .
\end{align*}
which in turn produces
\begin{align*}
&\|\partial_t(\widetilde\phi-\widetilde\phi_h)\|_{L^2(Q_T)}
+d_j^{-1}\|\nabla (\widetilde\phi-\widetilde\phi_h)\|_{L^2(Q_T)}
+d_j^{-1}\|\widetilde\phi-\widetilde \phi_h\|_{L^\infty((0,T);L^2)}
\nn\\
&\leq Cd_j^{-2}\| \widetilde\phi \|_{L^2(Q_T)}
+Cd_j^{-1}\| \nabla\widetilde\phi \|_{L^2(Q_T)}
+C\|\partial_t  \widetilde\phi \|_{L^2(Q_T)}
+Cd_j \| \nabla\partial_t\widetilde \phi \|_{L^2(Q_T)} .
\end{align*}
By applying Lemma \ref{lemlocEng0} to \refe{wdtzh2}-\refe{wdtzh3}
and using the inequality above,
we derive that
\begin{align*}
&\|\partial_t(\widetilde\phi_h-\phi_h)\|_{L^2(Q_j)}
+ d_j^{-1}\|\nabla(\widetilde \phi_h-\phi_h)\|_{L^2(Q_j)} \nn\\
&\leq  C_m\left(\|\nabla (\widetilde \phi_h-\phi_h)(0)\|_{L^2(\Omega_j'')}
+d_j^{-1}\|(\widetilde \phi_h-\phi_h)(0)\|_{L^2(\Omega_j'')} \right) \nn\\
&\quad
+ C_m\left(\frac{h}{d_j}\right)^{m}\big(
\|\partial_t(\widetilde \phi_h-\phi_h)\|_{L^2(Q_{j}'')}
+ d_j^{-1}\|\nabla(\widetilde \phi_h-\phi_h)\|_{L^2(Q_{j}'')}\big)\nn\\
&\quad
+ C_md_j^{-2} \| \widetilde \phi_h-\phi_h\|_{L^2(Q_{j}'')} \nn\\
&\leq  C_m\left(\|\nabla \phi_h(0)\|_{L^2(\Omega_j'')}
+d_j^{-1}\|\phi_h(0)\|_{L^2(\Omega_j'')} \right) \nn\\
&\quad
+ C_m\left(\frac{h}{d_j}\right)^{m}\big(
\|\partial_te\|_{L^2(Q_{j}'')}
+ d_j^{-1}\|\nabla e\|_{L^2(Q_{j}'')}\big)
+ C_md_j^{-2} \|e\|_{L^2(Q_{j}'')}  \nn\\
&\quad
+ C_m\left(\frac{h}{d_j}\right)^{m}\big(
\|\partial_t(\widetilde\phi-\widetilde\phi_h)\|_{L^2(Q_{j}'')}
+ d_j^{-1}\|\nabla (\widetilde\phi-\widetilde\phi_h)\|_{L^2(Q_{j}'')}\big) \nn\\
&\quad
+ C_md_j^{-1} \| \widetilde \phi-\widetilde\phi_h\|_{L^\infty L^2(Q_j'')} \nn\\
&\leq  C_m\left(\|\nabla \phi_h(0)\|_{L^2(\Omega_j'')}
+d_j^{-1}\|\phi_h(0)\|_{L^2(\Omega_j'')} \right) \nn\\
&\quad
+ C_m\left(\frac{h}{d_j}\right)^{m}\big(
\|\partial_te\|_{L^2(Q_{j}'')}
+ d_j^{-1}\|\nabla e\|_{L^2(Q_{j}'')}\big)
+ C_md_j^{-2} \|e\|_{L^2(Q_{j}'')} \nn\\
&\quad
+ C_md_j^{-2}\| \widetilde\phi \|_{L^2(Q_T)}
+C_md_j^{-1}\| \nabla\widetilde\phi \|_{L^2(Q_T)}
+C_m\|\partial_t  \widetilde\phi \|_{L^2(Q_T)} \nn\\
&\quad
+C_md_j \| \nabla\partial_t\widetilde \phi \|_{L^2(Q_T)}
\end{align*}
The last two inequalities imply
\begin{align}\label{etQjdj}
\|e_t\|_{L^2(Q_j)}
+ d_j^{-1}\|\nabla e\|_{L^2(Q_j)}
\leq
C_m\big(I_j(\phi_{h}(0))+X_j(\phi)
+H_j(e)+d_j^{-2}\|e\|_{L^2(Q_j''')}\big) .
\end{align}
We have proved that
any $e=\phi_h-\phi$ satisfying \refe{p32}
also satisfies \refe{etQjdj}.
Since
$e=\phi_h-\Pi_h\phi-(\phi-\Pi_h\phi)$
and $\phi(0)-\Pi_h\phi(0)=0$ in $\Omega_j'''$,
we can replace $\phi_h$ by $\phi_h-\Pi_h\phi$
and $\phi$ by $\phi-\Pi_h\phi$
in \refe{etQjdj}. Then \refe{3-1} follows immediately.
\,\endproof

\section{Proof of Lemma \ref{lemma2-1}}
\setcounter{equation}{0}

Now we turn back to the proof of Lemma \ref{lemma2-1}.

\subsection{
The proof of \refe{DxtGrEst}}\label{EGF}
In this subsection, we present several local energy estimates for the Green's
function, the regularized Green's function and the discrete Green's
function, which then are used to prove \refe{DxtGrEst}.
These energy estimates
will also be used to prove \refe{FFEst2} in the next subsection.
In this subsection we let $T=1$
and fix $x_0\in\Omega$. We write
$G$ and $\Gamma$ as abbreviations for
the functions $G(\cdot,\cdot,x_0)$ and
$\Gamma(\cdot,\cdot,x_0)$, respectively,
when there is no ambiguity.
We use the decomposition of Section \ref{SecGF}
for all $j\geq 1$ (not restricted to $j\leq J_*$)
and so we do not
require $h< R_0K_0^{-2}/(16C_*)$ in this subsection.

\begin{lemma}\label{GFEst1}
{\it
For the Green's functions $\Gamma, G$ and  $\Gamma_h$ defined in
\refe{GFdef}-\refe{GMhFdef}, we have the following estimates:
\begin{align}
&\sum_{l,k=0}^2d_j^{2l+k-1+N/2}
(\|\nabla^k\partial_{t}^lG(\cdot,\cdot,x_0)\|_{L^2(Q_j(x_0))}
+\|\nabla^k\partial_{t}^l\Gamma(\cdot,\cdot,x_0)\|_{L^2(Q_j(x_0))})
\leq C, \label{GFest01}\\
&\|\nabla^2G(\cdot,\cdot,x_0) \|_{L^{\infty,2}(\cup_{k\leq
j}Q_k(x_0))}\leq Cd_j^{-N/2-2} \label{GFest03}
,\\
&\|
\nabla_{x_0}\partial_tG(\cdot,\cdot, x_0 )\|_{L^{\infty}(Q_j(x_0))}\leq Cd_j^{-N-3},
\label{DxtGL1}\\
&\|\partial_tG(\cdot,\cdot,x_0)\|_{L^1(\Omega\times(1,\infty))}
+\|t\partial_{tt}G(\cdot,\cdot,x_0)\|_{L^1(\Omega\times(1,\infty))} \nn\\
&\quad
+\|\partial_t\Gamma(\cdot,\cdot,x_0)\|_{L^1(\Omega\times(1,\infty))}
+\|t\partial_{tt}\Gamma(\cdot,\cdot,x_0)\|_{L^1(\Omega\times(1,\infty))}
\leq C ,\label{GFest0423}
\end{align}
}
\end{lemma}

\noindent{\it Proof.}$\,\,$
For the given $x_0$ and $j$,  we define a coordinate transformation
$x-x_0=d_j\widetilde x$ and $t=d_j^2\widetilde t$, and
define $\widetilde
G(\widetilde t,\widetilde x):=G(t,x,x_0)$, $\widetilde
a(\widetilde x):=a(x)$.
Via the coordinates transformation, we assume that
the sets $Q_j$, $Q_j'$, $\Omega_j$,
$\Omega_j'$ and $\Omega$ are transformed to
$\widetilde Q_j$, $\widetilde Q_j'$, $\widetilde \Omega_j$,
$\widetilde \Omega_j'$ and $\widetilde \Omega$, respectively.
%
Let $0\leq \widetilde\omega_i(x,t)\leq 1$, $i=0,1,2,3$,
be smooth cut-off functions which vanishes outside
$\widetilde Q_j'$ and equals $1$ in $\widetilde Q_{j}$.
Moreover, $\widetilde\omega_{i}$ equals $1$
at the points where $\widetilde\omega_{i+1}\neq 0$,
and $|\nabla \widetilde\omega_i|\leq C$,
$|\partial_t\widetilde\omega_i|\le C$ for $i=0,1,2,3$.
Since $ \cup_{k\geq j}\widetilde\Omega_k'
\cup \widetilde\Omega_*$ is of unit size,
there exists a convex domain
$\widetilde D=B_{\rho}(z)\cap \widetilde\Omega
\supset \cup_{k\geq j}\widetilde\Omega_k'
\cup \widetilde\Omega_*$,
with $4\leq \rho\leq 8K_0^2$, which belongs to one
of the following cases
(there are only a finite
number of shapes for $\widetilde D$):

(i) $z\in \widetilde\Omega$, $\rho=4$, and
$B_4(z)$ has no intersection with the boundary of
$ \widetilde\Omega$, thus $B_\rho(z)\cap\widetilde\Omega=B_\rho(z)$,

(ii) $z$ is on a face of $ \widetilde\Omega$, $\rho=8$
and $B_8(z)$ has no intersection with other faces of $\widetilde\Omega$,
thus $B_\rho(z)\cap\Omega$ is a half ball,

(iii) $z$ is on an edge of $\widetilde\Omega$,
$\rho=8K_0$ and $B_{8K_0}(z)$
has no intersection with any closed faces of $\widetilde\Omega$
which do not contain this edge,
thus $B_\rho(z)\cap\Omega$ is the intersection of a ball with
a sector spanned by the edge,

(iv) $z$ is a corner of $ \widetilde\Omega$ and $\rho=8K_0^2<R_0$,
and $B_\rho(z)\cap \widetilde\Omega$ coincides with the intersection
of the ball $B_\rho(z)$ with the cone spanned by the corner $z$.

Note that $\widetilde D\times(0,16)$ contains $\widetilde Q_j'$,
and consider $\widetilde\omega_i\widetilde G$, $i=1,2$,
which are solutions of
\begin{align}\label{omegaG1}
\partial_{\widetilde t}(\widetilde\omega_1\widetilde G)
-\nabla_{\widetilde x} \cdot(\widetilde a
\nabla_{\widetilde x} (\widetilde\omega_1\widetilde G))
= \widetilde\omega_0\widetilde G\partial_{\widetilde t}\widetilde\omega_1
+ \widetilde\omega_0\widetilde G\nabla_{\widetilde x}\cdot\big(\widetilde
a\nabla_{\widetilde x} \widetilde\omega_1\big)
-\nabla_{\widetilde x} \cdot(2\widetilde
a \widetilde\omega_0\widetilde G\nabla_{\widetilde x} \widetilde\omega_1)
\end{align}
and
\begin{align}\label{omegaG2}
\partial_{\widetilde t}(\widetilde\omega_2\widetilde G)
-\nabla_{\widetilde x} \cdot(\widetilde a
\nabla_{\widetilde x} (\widetilde\omega_2\widetilde G))
=  \widetilde\omega_1\widetilde
G\partial_{\widetilde t}\widetilde\omega_2
+  \widetilde\omega_1\widetilde
G \nabla_{\widetilde x}\cdot\big(\widetilde
a\nabla_{\widetilde x} \widetilde\omega_2\big)
-\nabla_{\widetilde x} \cdot(2\widetilde
a \widetilde\omega_1\widetilde
G\nabla_{\widetilde x} \widetilde\omega_2)
\end{align}
in the domain $\widetilde D\times(0,16)$, respectively,
both with zero boundary/initial conditions.
Since $\widetilde D$ is a convex domain,
for $p=N+\alpha$ and $p_0=Np/(N+p)$ so that
$W^{1,p}(\widetilde D)\hookrightarrow L^\infty(\widetilde D)$ and
$W^{1,p_0}(\widetilde D)\hookrightarrow L^p(\widetilde D)$,
the standard $L^p((0,16);W^{1,p}(\widetilde D))$ estimate
of \refe{omegaG1} (the inequality \refe{LpW2quf}-\refe{LpW1quf} with $p=q$,
Lemma \ref{LemMaxLp}) gives
\begin{align*}
&\|\nabla_{\widetilde x}(\widetilde\omega_1
\widetilde G)\|_{L^{p}((0,16);L^p(\widetilde D))} \\
& \leq
C\|\widetilde\omega_0\widetilde G \|_{L^{p}((0,16);L^{p_0}(\widetilde D))}
+C\|\widetilde\omega_0\widetilde G \nabla_{\widetilde x}\widetilde a \|_{L^{p}((0,16);L^{p_0}(\widetilde D))}
+C\|\widetilde\omega_0\widetilde G \|_{L^{p}((0,16);L^p(\widetilde D))} \\
&\leq
C\|\widetilde\omega_0\widetilde G \|_{L^{p}((0,16);L^p(\widetilde D))}
(C+C\|\nabla_{\widetilde x}\widetilde a \|_{L^N(\widetilde D)}) \\
&\leq
C\|\widetilde\omega_0\widetilde G \|_{L^{p}((0,16);L^p(\widetilde D))}
,
\end{align*}
and the maximal $L^p$ regularity of \refe{omegaG2} yields that
(see inequality \refe{LpW2quf},
Lemma \ref{LemMaxLp})
\begin{align*}
&\|\partial_{\widetilde t}(\widetilde\omega_2\widetilde G ))
\|_{L^{p}((0,16);L^p(\widetilde D))} +\|\nabla_{\widetilde x} \cdot(\widetilde a
\nabla_{\widetilde x} (\widetilde\omega_2\widetilde G ))
\|_{L^{p}((0,16);L^p(\widetilde D))}  \\
&
\leq
C\|\widetilde \omega_1\widetilde G\|_{L^{p}((0,16);L^p(\widetilde D))}
+C\|\widetilde \omega_1\widetilde G\nabla_{\widetilde x}\widetilde a\|_{L^{p}((0,16);L^p(\widetilde D))} \\
&\quad
+C\|\nabla_{\widetilde x}(\widetilde \omega_1\widetilde G  )\|_{L^{p}((0,16);L^p(\widetilde D))}  \\
&\leq
C\|\widetilde \omega_1\widetilde G\|_{L^{p}((0,16);L^{\infty}(\widetilde D))}
(C+C\|\nabla_{\widetilde x}\widetilde a\|_{L^{p}(\widetilde D)}  )
+C\|\nabla_{\widetilde x}(\widetilde \omega_1\widetilde G  )\|_{L^{p}((0,16);L^p(\widetilde D))} \\
&\leq
C\|\nabla_{\widetilde x}(\widetilde \omega_1\widetilde G  )\|_{L^{p}((0,16);L^p(\widetilde D))}  .
\end{align*}
By using \refe{H2Reg0}-\refe{W1inftyReg}, we have
\begin{align*}
&\|\nabla_{\widetilde x}(\widetilde\omega_2\widetilde G )
\|_{L^{\infty}(\widetilde D)}
+\sum_{k=0}^2
\|\nabla_{\widetilde x}^k(\widetilde\omega_2\widetilde G)\|_{L^{2}(\widetilde D)}\leq
C\|\nabla_{\widetilde x} \cdot(\widetilde a
\nabla_{\widetilde x} (\widetilde\omega_2\widetilde G ))
\|_{L^{p}(\widetilde D)} .
 \end{align*}
The last three inequalities imply that
\begin{align*}
&\|\nabla_{\widetilde x}(\widetilde\omega_2\widetilde G )
\|_{L^2((0,16);L^{\infty}(\widetilde D))}
+\sum_{k=0}^2\|\nabla_{\widetilde x}^k
(\widetilde\omega_2\widetilde G)\|_{L^2((0,16);L^2(\widetilde D))} \nn\\
&\quad
+\|\partial_{\widetilde t}(\widetilde\omega_2\widetilde G )
\|_{L^{N+\alpha}((0,16);L^{N+\alpha}(\widetilde D))} \\
&\leq
C\|\widetilde\omega_0\widetilde G \|_{L^{N+\alpha}((0,16);L^{N+\alpha}(\widetilde D))} .
 \end{align*}

Similarly,
replacing $\widetilde G$ by $\partial_{\widetilde t}\widetilde G$
and $\partial_{\widetilde t}^2\widetilde G$
in the above estimates, respectively, one can derive that
\begin{align*}
&\sum_{l=0}^2\|\nabla_{\widetilde x}(\widetilde\omega_3
\partial_{\widetilde t}^l\widetilde G )
\|_{L^2((0,16);L^{\infty}(\widetilde D))}
+\sum_{l,k=0}^2\|\nabla_{\widetilde x}^k(
\widetilde\omega_3\partial_{\widetilde t}^l
\widetilde G)\|_{L^{2}((0,16);L^2(\widetilde D))}  \\
&\leq
C\|\widetilde\omega_0\widetilde G \|_{L^{N+\alpha}((0,16);L^{N+\alpha}(\widetilde D))}
\leq
C\|\widetilde\omega_0\widetilde G \|_{L^{\infty}((0,16);L^\infty(\widetilde D))} .
\end{align*}
Since $\widetilde\omega_3\widetilde G\equiv 0$ at $t=0$,
it follows that
\begin{align*}
\|\nabla_{\widetilde x}\partial_{\widetilde t}(\widetilde\omega_3
\widetilde G)\|_{L^\infty((0,16);L^\infty(\widetilde D))}
&\leq C
\|\nabla_{\widetilde x}\partial_{\widetilde t}^2
(\widetilde\omega_3\widetilde G)\|_{L^2((0,16);L^\infty(\widetilde D))}\\
&
\leq C\|\widetilde\omega_0\widetilde G \|_{L^{\infty}((0,16);L^\infty(\widetilde D)))} , 
\end{align*}
and
\begin{align*}
\sum_{k=0}^2\|\nabla_{\widetilde x}^k
(\widetilde\omega_3\widetilde G)\|_{L^\infty((0,16);L^2(\widetilde D))}
&\leq C\sum_{k=0}^2\|\partial_{\widetilde t} \nabla_{\widetilde x}^k
(\widetilde\omega_3\widetilde G)\|_{L^2((0,16);L^2(\widetilde D))} \\
&\leq C\|\widetilde\omega_0\widetilde G \|_{L^{\infty}((0,16);L^\infty(\widetilde D))}
\, .
\end{align*}
Moreover, from the last three inequalities, we have
\begin{align*}
&\|\nabla_{\widetilde x} \partial_{\widetilde t}
\widetilde G \|_{L^{\infty}(\widetilde Q_j)}
+\sum_{k=0}^2\|\nabla_{\widetilde x}^k
\widetilde G\|_{L^{\infty,2}(\widetilde Q_j)}
+\sum_{l,k=0}^2\|\partial_{\widetilde t}^l
\nabla_{\widetilde x}^k\widetilde G\|_{L^2(\widetilde Q_j)}
\leq  C\|\widetilde G\|_{L^{\infty}(\widetilde Q_j')}  .
\end{align*}

Transforming back to the $(x,t)$ coordinates,
we see from the last two inequalities that
\begin{align*}
&d_j^{3}\|\nabla\partial_t G \|_{L^{\infty}(Q_j)}
+\sum_{k=0}^2 d_j^{k-N/2}\| \nabla^k G\|_{L^{\infty,2}(Q_j)}
+\sum_{l,k=0}^2 d_j^{2l+k-1-N/2}\|\partial_t^l\nabla^k G\|_{L^2(Q_j)} \\
&
 \leq C\| G \|_{L^{\infty}(Q_j')}\leq Cd_j^{-N} ,
\end{align*}
where we have used \refe{FEstP} in the last inequality.
By the symmetry of $G$ with respect to $x$ and $x_0$
we also get
$$
d_j^{3}\|\nabla_{x_0}\partial_t G (\cdot,\cdot,x_0)\|_{L^{\infty}(Q_j)}
\leq C\| G \|_{L^{\infty}(Q_j')}\leq Cd_j^{-N} .
$$
This proves \refe{GFest01}-\refe{DxtGL1}
for $G$.

By using the expression
\begin{align}
\Gamma(x,t\, ;x_0)
=\int_\Omega G(x,t\, ;y)\widetilde\delta(y,x_0)\d y .
\end{align}
one can derive the same estimates for $\Gamma$:
\begin{align*}
\|\nabla^k\partial_t^l\Gamma(\cdot,\cdot \, ; x_0)\|_{L^2(Q_j)}
&\leq\int_\Omega \|\nabla^k\partial_t^lG(\cdot,\cdot \, ;y)\|_{L^2(Q_j)}
|\widetilde\delta(y,x_0)|\d y \\
&\leq \int_\Omega Cd_j^{-2l-k+1-N/2}|\widetilde\delta(y,x_0)|\d y
\leq Cd_j^{-2l-k+1-N/2} .
\end{align*}

Finally, we note that
\begin{align} \label{Gammatt}
&\|\partial_t\Gamma (t,\cdot,x_0)\|_{L^\infty}
+t\|\partial_{tt}\Gamma  (t,\cdot,x_0)\|_{L^\infty} \nn \\
&\leq C\|\partial_t\Gamma (t,\cdot,x_0)\|_{H^2}
+Ct\|\partial_{tt}\Gamma  (t,\cdot,x_0)\|_{H^2} \nn \\
&\leq C\|\partial_tA\Gamma (t,\cdot,x_0)\|_{L^2}
+Ct\|\partial_{tt}A \Gamma  (t,\cdot,x_0)\|_{L^2}
&&\mbox{[$H^2$ estimate, Lemma \ref{LemMaxLp}]} \nn \\
&= C\|\partial_{tt}\Gamma (t,\cdot,x_0)\|_{L^2}
+Ct\|\partial_{ttt}\Gamma (t,\cdot,x_0)\|_{L^2} \nn \\
&\leq Ct^{-2}\|\Gamma (t/2,\cdot,x_0)\|_{L^2}
&&\mbox{[semigroup estimate]} \nn  \\
&\leq Ct^{-2-N/2} ,
\end{align}
which implies the first part of (\ref{GFest0423}) and
the second part of (\ref{GFest0423})
(the estimates of $G$) can be proved in the same way.

The proof of Lemma \ref{GFEst1} is completed.
\,\endproof \medskip

When $\max(t^{1/2},|x-x_0|)<d_1$ the inequality
\refe{DxtGrEst} follows from \refe{DxtGL1}.
When $\max(t^{1/2},|x-x_0|)\geq d_1$,
the estimate
$\|\nabla_{x_0}\partial_tG(\cdot,\cdot, x_0 )\|_{L^{\infty}(Q_j(x_0))}\leq C$
can be proved directly (without using the scale transformation)
in the same way as above.

\subsection{Proof of \refe{FFEst2}}
\label{dka7}
The proof is also based on Lemma \ref{GFEst1}.

First we consider the case $h < h_0:=(R_0K_0^{-2}/(16C_*)$ and let $T=1$.
The basic energy estimates of the equations
\refe{GMFdef}-\refe{GMhFdef} yield
\begin{align*}
\|\partial_t\Gamma_h\|_{L^2(Q_T)}
+\|\partial_t\Gamma\|_{L^2(Q_T)}
&\leq \|\nabla\Gamma_h(0)\|_{L^2(\Omega)}
+\|\nabla\Gamma(0)\|_{L^2(\Omega)} \\
&
=\|\nabla P_h\widetilde\delta_{x_0}\|_{L^2(\Omega)}
+\|\nabla \widetilde\delta_{x_0}\|_{L^2(\Omega)}
\leq  Ch^{-1-N/2}
\end{align*}
and
\begin{align*}
\|\partial_{tt}\Gamma_h\|_{L^2(Q_T)}
+\|\partial_{tt}\Gamma\|_{L^2(Q_T)}
&\leq \|\nabla\partial_t\Gamma_h(0)\|_{L^2(\Omega)}
+\|\nabla\partial_t\Gamma(0)\|_{L^2(\Omega)}  \\
&=\|\nabla A_hP_h\widetilde\delta_{x_0}\|_{L^2(\Omega)}
+\|\nabla A\widetilde\delta_{x_0}\|_{L^2(\Omega)}
\leq Ch^{-3-N/2} ,
\end{align*}
which imply
\begin{align*}
\|\partial_tF\|_{L^2(Q_*)}  +\|t\partial_{tt}F\|_{L^2(Q_*)}
&\leq Ch^{-1-N/2}+Cd_{J_*}^2h^{-3-N/2} \\
&\leq Ch^{-1-N/2}+CC_*^2h^{-1-N/2} .
\end{align*}
Hence we have
\begin{align}\label{Bd31K}
&\|\partial_tF\|_{L^1( Q_T)}+\|t\partial_{tt}F\|_{L^1( Q_T)} \nn\\
&\leq Cd_{J_*}^{1+N/2}\big(\|\partial_tF\|_{L^2(Q_*)}  +\|t\partial_{tt}F\|_{L^2(Q_*)}\big) \nn\\
&\quad +\sum_{j}Cd_j^{1+N/2}
\big(\|\partial_tF\|_{L^2(Q_j)}  +\|t\partial_{tt}F\|_{L^2(Q_j)}\big) \nn\\
&\leq CC_*^{3+N/2}
+\sum_{j}Cd_j^{1+N/2}\big(\|\partial_tF\|_{L^2(Q_j)}  +\|t\partial_{tt}F\|_{L^2(Q_j)}\big) \nn\\
&\leq CC_*^{3+N/2}+C{\mathcal K} ,
\end{align}
where
\begin{align}\label{KdKj}
{\mathcal K}:=\sum_{j}d_j^{1+N/2}(
d_{j}^{-1}\|\nabla F\|_{L^2(Q_j)}
+\|\partial_tF\|_{L^2(Q_j)}
+d_{j}^{2}\|\partial_{tt}F\|_{L^2(Q_j)}) .
\end{align}

We proceed to estimate ${\mathcal K}$.
We set $e=F$ ($\phi_h=\Gamma_h$ and $\phi=\Gamma$)
and $e=\partial_tF$ ($\phi_h=\partial_t\Gamma_h$ and $\phi=\partial_t\Gamma$)
in \refe{p32} (Proposition \ref{LocEEst}), respectively,
and note that $\Gamma(0)=\partial_t\Gamma(0)=0$ on $\Omega_j'$.
We obtain that
\begin{align}
&d_{j}^{-1}\|\nabla F\|_{L^2(Q_j)}+\|\partial_tF\|_{L^2(Q_j)}
+d_{j}^{2}\|\partial_{tt}F\|_{L^2(Q_j)} \nn\\
&\leq
C(\widehat{I_j}+\widehat{X_j}+\widehat{H_j}+ d^{-2}_j\| F \|_{L^2(Q_j')})
\end{align}
where
\begin{align*}
&\widehat{I_j}=\|\nabla P_h\widetilde\delta_{x_0}\|_{L^2(\Omega_j')}
+d_j^{-1}\|P_h\widetilde\delta_{x_0}\|_{L^2(\Omega_j')} \\
&\qquad 
+d_{j}^{2}\|\nabla A_hP_h\widetilde\delta_{x_0}\|_{L^2(\Omega_j')}
+d_j \|A_hP_h\widetilde\delta_{x_0}\|_{L^2(\Omega_j')},
\\
&\widehat{X_j}=
d_j\|\nabla\partial_t(\Pi_h\Gamma-\Gamma) \|_{L^2(Q_j')}
+\|\partial_t(\Pi_h\Gamma-\Gamma)\|_{L^2(Q_j')}
+d_j^{-1}\|\nabla(\Pi_h\Gamma-\Gamma)\|_{L^2(Q_j')}\\
&\qquad
+ d_j^{-2}\|\Pi_h\Gamma-\Gamma\|_{L^2(Q_j')} +
d_j^3\| \nabla\partial_{tt}(\Pi_h\Gamma-\Gamma) \|_{L^2(Q_j')}
+d_j^2\|\partial_{tt}(\Pi_h\Gamma-\Gamma)\|_{L^2(Q_j')} \\
&\widehat{H_j}=\big(h/d_j\big)^m\big(  \|\partial_{t}F\|_{L^2(Q_j')}
+d_j^{-1}\|\nabla F\|_{L^2(Q_j')}
+d_j^{2}\|\partial_{tt}F \|_{L^2(Q_j')}
+d_j \|\nabla\partial_{t}F\|_{L^2(Q_j')} \big)
\, .
\end{align*}
By noting the exponential decay of
$ P_h\widetilde \delta_{x_0}(y) $ (see (P2) in section \ref{Sec2-2})
we derive
\begin{align*}
&\widehat{I_j} \leq Chd_{j}^{-2-N/2},\\
&\widehat{X_j}
\leq (d_j h+ h^2)\|\nabla^2\partial_{t}\Gamma\|_{L^2(Q_j'')}
+(d_j^{-1}h+ d_j^{-2}h^2)\|\nabla^2\Gamma\|_{L^2(Q_j'')} \\
&\quad 
+( d_j^3 h+ d_j^{ 2}h^2)\|\nabla^2\partial_{tt}\Gamma\|_{L^2(Q_j'')} \\
&~
\leq C hd_j^{-2-N/2} ,
\quad 
\mbox{[by using Lemma \ref{GFEst1}]} \\
&\widehat{H_j}
\leq \big(h/d_j\big)^m \big(  \|\partial_{t}F\|_{L^2(Q_T)}
+d_j^{-1}\|\nabla F\|_{L^2(Q_T)}
+d_j^{2}\|\partial_{tt}F \|_{L^2(Q_T)}
+d_j \|\nabla\partial_{t}F\|_{L^2(Q_T)} \big) \\
&~
\leq C\big(h/d_j\big)^m  \big(\|P_h\widetilde\delta_{x_0} \|_{H^1}
+d_j^{-1}\|P_h\widetilde\delta_{x_0} \|_{L^2}
+d_j^{2}\|A_hP_h\widetilde\delta_{x_0} \|_{H^1}
+d_j \|A_hP_h\widetilde\delta_{x_0} \|_{L^2}\big)
\\
&~
\leq C\big(h/d_j\big)^m \big(h^{-1-N/2}
+d_j^{-1}h^{-N/2}+d_j^{2}h^{-3-N/2}
+d_j h^{-2-N/2}\big)
\quad  \mbox{[by (P3)-(P4)]}\\
&~
\leq Chd_j^{-2-N/2},\quad \mbox{[by choosing $m=4+N/2$]}
\end{align*}
Therefore, by \refe{KdKj},
\begin{align}\label{dlj6}
{\mathcal K}
&\leq \sum_{j} Chd_j^{-1}+\sum_{j}Cd_j^{-1+N/2}\|F\|_{L^2(Q_j')} \nn\\
&\leq C+C\sum_{j}d_j^{-1+N/2}\|F\|_{L^2(Q_j')} .
\end{align}

To estimate $\|F\|_{L^2(Q_j')}$, we apply a duality argument.
Let $w$ be the solution of the backward parabolic equation
\begin{align*}
\left\{\begin{array}{ll}
-\partial_tw+Aw=v &\mbox{in}\,\,\,\Omega ,\\
w=0 &\mbox{on}\,\,\,\partial\Omega,\\
w(T)=0 &\mbox{in}\,\,\,\Omega ,
\end{array}\right.
\end{align*}
where $v$ is a function which is supported on $Q_j'$ and
$\|v\|_{L^2(Q_j')}=1$. Multiplying the above equation by $F$,
with integration by parts we get
\begin{align}\label{dka6}
\iint_{ Q_T} Fv dx dt
=(F(0),w(0)) + \iint_{ Q_T} \partial_tF w dx dt
+\sum_{i,j=1}^N \iint_{ Q_T} a_{ij}\partial_j F \partial_i w dx dt ,
\end{align}
where
\begin{align*}
(F(0),w(0))&=(P_h\widetilde\delta_{x_0}-\widetilde\delta_{x_0},w(0))\\
&=(P_h\widetilde\delta_{x_0}-\widetilde\delta_{x_0},w(0)-I_hw(0))\\
&=
(P_h\widetilde\delta_{x_0},w(0)-I_hw(0))_{\Omega_j''}
+(P_h\widetilde\delta_{x_0}-\widetilde\delta_{x_0},
w(0)-I_hw(0))_{(\Omega_j'')^c}\\
&:={\mathcal I}_1+{\mathcal I}_2 .
\end{align*}
By using the exponential decay of $P_h\widetilde\delta_{x_0}$
(see (P4) of section 2) and the local approximation
property (see (P2) of section 2), we derive that
\begin{align}
&|{\mathcal I}_1|\leq
Ch\|P_h\widetilde\delta_{x_0}\|_{L^2(\Omega_j'')}
(d_j^{-1}\|w(0)\|_{L^2(\Omega)}  + \|\nabla w(0)\|_{L^2(\Omega)}  )
\nn\\
&\quad\,\,\, \leq
Ch^{-N/2+1}e^{-Cd_j/h} (d_j^{-1}\|v\|_{L^{2(N+2)/(N+4)}(Q_j')}  +\|v\|_{L^2(Q_j')})
\nn\\
&\quad\,\,\, \leq
Ch^{-N/2+1}e^{-Cd_j/h}\|v\|_{L^2(Q_j')}
\nn\\
&\quad\,\,\, \leq
C(d_j/h)^{1+N/2}e^{-Cd_j/h}h^{2}d_j^{-1-N/2}
\nn\\
&\quad\,\,\, \leq
Ch^{2}d_j^{-1-N/2 }, \label{SF2}\\[5pt]
&|{\mathcal I}_2|\leq C\| \widetilde\delta_{x_0}\|_{L^{2} }
\|w(0)-I_hw(0)\|_{L^{2}((\Omega_j'')^c)}
\nn\\
&\quad\,\,\, \leq
Ch^{2-N/2}\sum_{k=0}^2 d_j^{k-2}\|\nabla^k w(0)\|_{L^{2}((\Omega_j'')^c)} .
\label{SF22}
\end{align}

To estimate $\|\nabla^kw(0)\|_{L^{2}((\Omega_j'')^c)}$,
we let $W_j$ be a set containing $(\Omega_j'')^c$ but its
distance to $\Omega_j'$ is larger than $d_j/8$. Since
$$
\nabla^k_xw(x,0)
=\int_0^{T}\int_{\Omega}
\nabla^k_xG(s,x,y)v(y,s)\d y\d s  ,
$$
by noting the fact
$$
|x-y|+s^{1/2}\geq  d_j /8
\quad\mbox{for $x\in W_j$ and $(y,s)\in Q_j'$}
$$
and using (\ref{GFest03}), we
further derive
\begin{align}
&\sum_{k=0}^2d_j^{k-2}\|\nabla^kw(\cdot,0)\|_{L^{2}(W_j)} \nn\\
&\leq C\sum_{k=0}^2d_j^{k-2} \sup_{y\in\Omega}\|\nabla^k
G(\cdot,\cdot,y)\|_{L^{\infty,2}(\bigcup_{k\leq
j-3}Q_k(y))}\|v\|_{L^{1}(Q_j')}\nn\\
&\leq \sum_{k=0}^2d_j^{k-2} d_j^{-N -k+N/2}\|v\|_{L^{1}(Q_j')}
 \nn\\
&\leq \sum_{k=0}^2d_j^{k-2} d_j^{-N -k+N/2}d_j^{N/2+1}\|v\|_{L^2(Q_j')}
=C d_j^{ -1 } . \label{SF3}
\end{align}
From (\ref{SF2})-(\ref{SF3}), we see that the
first term on the right-hand side of \refe{dka6} is bounded by
\begin{align}
|(F(0),w(0))| \leq  Ch^{2}d_j^{-N/2 -1}+Ch^{2-N/2}d_j^{ -1}
 \leq  Ch^{2-N/2}d_j^{-1} ,
\label{f69}
\end{align}
and the rest terms are bounded by
\begin{align}\label{sd80}
\iint_{ Q_T} \partial_tF w dx dt
& +\sum_{i,j=1}^N \iint_{ Q_T} a_{ij}\partial_j F \partial_iw dx dt
\nn \\
&= \iint_{ Q_T} \partial_tF (w-\Pi_hw) dx dt
+\sum_{i,j=1}^N \iint_{ Q_T} a_{ij}\partial_j F \partial_i (w-\Pi_hw)dx dt
\nn\\
&\leq
\sum_{*,i} C\|\nabla^2w\|_{L^2(Q_i')} (h^2\|\partial_tF\|_{L^2(Q_i)}
+ h\|\nabla F\|_{L^2(Q_i)} )
 .
\end{align}

Moreover, to estimate $\|\nabla^2w\|_{L^2(Q_i')}$
we consider the expression
$$
\nabla^2_x w(x,t)
=\int_0^{T}\int_{\Omega}
\nabla^2_x G(s-t,x,y)v(y,s)1_{s>t}\,\d y\d s .
$$
For $i\leq j-3$ (so that $d_i>d_j$), we see that $w(x,t)=0$ for
$t>16d_j^2$ (because $v$ is supported in $Q_j'$); $d_i/2\leq |x-y|\leq 4 d_i$
and $s-t<d_i^2$ for $t< 16d_j^2$, $(x,t)\in Q_i$ and
$(y,s)\in Q_j'$.
Therefore, $(x,t)\in Q_i'(y)$ and we obtain
\begin{align*}
\|\nabla^2 w\|_{L^2(Q_i')}
&\leq \sup_{y}\|\nabla^2G(\cdot,\cdot,y)\|_{L^2(Q_i'(y))}
\|v\|_{L^1(Q_j')}\\
&\leq
Cd_i^{-N/2-1}d_j^{N/2+1}\|v\|_{L^2(Q_i')}\\
&\leq
C(d_j/d_i)^{N/2+1}
\leq \frac{Cd_j}{d_i} .
\end{align*}
For $i\geq j+3$ (so that $d_i< d_j$),
$\max(|s-t|^{1/2},|x-y|)\geq d_{j+1}$ for $(x,t)\in Q_i$ and therefore,
\begin{align*}
\|\nabla^2w\|_{L^2(Q_i')}
&\leq \ \sup_{y\in \Omega}\|\nabla^2G(\cdot,\cdot,y)
1_{\bigcup_{k\leq j+1}Q_{k}(y)}\|_{L^2(Q_i')}
\|v\|_{L^1(Q_j')}\\
&\leq Cd_i \sup_{y}\|\nabla^2G(\cdot,\cdot,y)
\|_{L^{\infty,2}(\bigcup_{k\leq j+1}Q_{k}(y))}
\|v\|_{L^2(Q_j')}d_j^{N/2+1} \\
&\leq C d_i d_j^{-N -2+N/2}d_j^{N/2+1}
=\frac{ C d_i}{d_j}  .
\end{align*}
Finally for $|i-j|\leq 2$, applying the standard energy estimate leads to
$$\|w\|_{L^2((0,T);H^2)}\leq C\|v\|_{L^2(Q_T)}=C.$$
Combining the three cases, we have
\begin{align}\label{sd808}
\|\nabla^2w\|_{L^2(Q_i')}\leq C \min \big(d_i/d_j,d_j/d_i\big) :=C m_{ij}.
\end{align}
Substituting \refe{f69}-\refe{sd808} into (\ref{dka6}) gives the estimate
\begin{align}
\|F\|_{L^2(Q_j')}\leq Ch^2d_j^{-N/2-1}
 +C\sum_{*,i}m_{ij}
(h^2\|\partial_tF\|_{L^2(Q_i)}
+h\|\nabla F\|_{L^2(Q_i)}) ,
\end{align}
which together with \refe{dlj6} implies
\begin{align*}
{\mathcal K}&\leq C+C\sum_j \left(\frac{h}{d_j}\right)^{2-N/2} +C\sum_j d_j^{N/2-1}
  \sum_{*,i}m_{ij}\big(h^2\|\partial_tF\|_{L^2(Q_i)}
+h\|\nabla F\|_{L^2(Q_i)}\big)\\
&\leq
C  +C\sum_{*,i}\left(h^2\|\partial_tF\|_{L^2(Q_i)}
+h\|\nabla F\|_{L^2(Q_i)} \right)\sum_jd_j^{N/2-1}m_{ij}\\
&\leq
C +C\sum_{*,i}\left(h^2\|\partial_tF\|_{L^2(Q_i)}
+h\|\nabla F\|_{L^2(Q_i)}\right) d_i^{N/2-1}\\
&\leq C +C\left(h^2\|\partial_tF\|_{L^2(Q_*)}
+h\|\nabla F\|_{L^2(Q_*)}\right)(C_* h )^{N/2-1} \\
&\quad +C\sum_id^{1+N/2}_i
\left(\|\partial_tF\|_{L^2(Q_i)}+d_i^{-1}\|\nabla F\|_{L^2(Q_i)}\right
)\left(\frac{h}{d_i}\right) \\
&\leq
C +CC_*^{-1 +N/2 }+C\sum_id^{1+N/2}_i\left(\|\partial_tF\|_{L^2(Q_i)}+d_i^{-1}\|\nabla F\|_{L^2(Q_i)}\right
)\left(\frac{h}{d_i}\right) \\
&\leq  C_2+C_2 C^{-1 +N/2 }_* +C_2C^{-1}_*{\mathcal K}
\end{align*}
for some positive constant $C_2$.
By choosing
\begin{align}\label{Cstar}
C_*= \max(10,10\kappa,R_0K_0^{-2}/8) +2C_2,
\end{align}
the above inequality shows that ${\mathcal K} \le C$.

Returning to \refe{Bd31K}, the boundedness of ${\mathcal K}$
implies
\begin{align}\label{Bd31K2}
&\|\partial_tF\|_{L^1( Q_T)}+\|t\partial_{tt}F\|_{L^1( Q_T)}\leq C.
\end{align}
From \refe{KdKj} we also see that, the boundedness
of ${\mathcal K}$ implies
\begin{align*}
\|\partial_tF\|_{L^2( Q_j)}
+\|t\partial_{tt}F\|_{L^2(Q_j)}\leq  Cd_j^{-1-N/2} .
\end{align*}
Since $\Omega\times(1/4,1)\subset \cup_{d_j\geq 1/2}Q_j$,
it follows that
\begin{align*}
&\|\partial_tF\|_{L^2(\Omega\times(1/4,1))}^2
+\|t\partial_{tt}F\|_{L^2(\Omega\times(1/4,1))}^2 \\
&\leq \sum_{d_j\geq 1/2}\big(\|\partial_tF\|_{L^2(Q_j)}^2
+\|t\partial_{tt}F\|_{L^2(Q_j)}^2\big)\\
&
\leq  \sum_{d_j\geq 1/2} C  d_j^{-2-N}  \leq C .
\end{align*}
The above inequality and (\ref{Gammatt}) imply
$$
\|\partial_t\Gamma_h\|_{L^2(\Omega\times(1/4,1))}
+\|t\partial_{tt}\Gamma_h\|_{L^2(\Omega\times(1/4,1))}  \leq C .
$$

Furthermore,  differentiating the equation (\ref{GMhFdef}) with respect to
$t$ and multiplying the result by $\partial_t\Gamma_h$ give
\begin{align*}
&\frac{\d}{\d t}\|\partial_t\Gamma_h(t,\cdot,x_0)\|_{L^2}^2
+c_0\|\partial_t\Gamma_h(t,\cdot,x_0)\|_{L^2}^2  \\
&\leq \frac{\d}{\d t}\|\partial_t\Gamma_h(t,\cdot,x_0)\|_{L^2}^2
+(A_h\partial_t\Gamma_h(t,\cdot,x_0),\partial_t\Gamma_h(t,\cdot,x_0))
= 0 ,\quad\mbox{for}\quad t\geq 1,
\end{align*}
which further shows that
$$ \|\partial_t\Gamma_h(t,\cdot,x_0)\|_{L^2}^2\leq
Ce^{-c_0(t-1)}\|\partial_t\Gamma_h(\cdot,\cdot,x_0)\|_{L^2((1/4,1);L^2(\Omega))}^2
\leq Ce^{-c_0t} \quad\mbox{for}\,\,\, t\geq 1. $$
In a similar way one can derive
$\|\partial_{tt}\Gamma_h(t,\cdot,x_0)\|_{L^2} \leq Ce^{-c_0t}$
for $t\geq 1$.
These inequalities together with \refe{GFest0423} imply
\begin{align}
\|\partial_tF (\cdot,\cdot,x_0)\|_{L^1(\Omega\times(1,\infty))}
+\|t\partial_{tt}F (\cdot,\cdot,x_0)\|_{L^1(\Omega\times(1,\infty))}\leq C ,
\end{align}
which together with \refe{Bd31K2} leads to
(\ref{FFEst2}) for the case
$h<h_0:=R_0K_0^{-2}/(16C_*)$ with $C_*$ being given by \refe{Cstar}.

Secondly when $h\geq h_0$, the decomposition in subsection 2.3 is not needed and the energy estimates of
\refe{GMFdef}-\refe{GMhFdef} yield
\begin{align*}
&\|\partial_t\Gamma_h\|_{L^2(Q_T)}
+\|\partial_t\Gamma\|_{L^2(Q_T)}
\leq \|\nabla\Gamma_h(0)\|_{L^2(\Omega)}
+\|\nabla\Gamma(0)\|_{L^2(\Omega)} \leq  Ch_0^{-1-N/2} \\
&\|\partial_{tt}\Gamma_h\|_{L^2(Q_T)}
+\|\partial_{tt}\Gamma\|_{L^2(Q_T)}
\leq \|\nabla\partial_t\Gamma_h(0)\|_{L^2(\Omega)}
+\|\nabla\partial_t\Gamma(0)\|_{L^2(\Omega)} \leq Ch_0^{-3-N/2} ,
\end{align*}
which imply
\begin{align}
&\int_0^1\int_\Omega\big( |\partial_tF(t,x,x_0) \big|
+\big|t\partial_{tt}F(t,x,x_0 ) |\big)\d x\d t \leq C .
\end{align}
Since both $\|\partial_{t}\Gamma_h(t)\|_{L^2(\Omega)}
+\|\partial_{tt}\Gamma_h(t)\|_{L^2(\Omega)}$
and $\|\partial_{t}\Gamma(t)\|_{L^2(\Omega)}
+\|\partial_{tt}\Gamma(t)\|_{L^2(\Omega)}$  decay
exponentially as $t\rightarrow \infty$,
it follows that \refe{FFEst2} still holds
when $h\geq h_0$.

The proof of Lemma \ref{lemma2-1} is completed.
~\endproof

\section{Conclusion}
In this paper we have proved that the
discrete elliptic operator $-A_h$
generates a bounded analytic semigroup
and has the maximal $L^p$ regularity, uniformly with respect to $h$,
in arbitrary convex polygons and polyhedra under the
regularity assumption $a_{ij}\in W^{1,N+\alpha}$.
We have assumed the quasi-uniformity of the triangulation, and
analysis of the problem under non-quasi-uniform triangulations
remains open. As far as we know,
only the analytic semigroup estimate \refe{STLEst}
and its equivalent resolvent estimate
were studied with an extra logarithmic factor
for some special cases of non-quasi-uniform triangulations, see \cite{CT01, Thomee07}.
The discrete maximal regularity estimates
(\ref{LpqSt3}) and (\ref{LpqSt2}) have not been established
with more general triangulations even
in smooth settings.

\section*{Appendix: The proof of Lemma 2.1}
\renewcommand{\thelemma}{A.\arabic{lemma}}
\renewcommand{\thetheorem}{A.\arabic{theorem}}
\renewcommand{\theequation}{A.\arabic{equation}}
\setcounter{equation}{0}

{\it Proof of \refe{H2Reg0}:}$\,\,$
The inequality \refe{H2Reg0} is similar to
Theorem 3.1.3.1 of \cite{Grisvard},
which was proved by using the local energy inequality of
Lemma 3.1.3.2, and the lemma was proved under the assumption
$a_{ij}\in W^{1,\infty}(\Omega)$, where $\Omega$ is a convex domain.
In the following, we show that
this assumption can be relaxed to
$a_{ij}\in W^{1,N+\alpha}(\Omega)
\hookrightarrow C^{\gamma}(\overline\Omega)$,
where $\gamma=\alpha/(N+\alpha)$

\begin{lemma}\label{LemLocH2}
{\it
If $\Omega$ is convex and $a_{ij}\in W^{1,N+\alpha}(\Omega)$, then
each point $y\in \Omega$ has a neighborhood
$B_R(y)\cap\Omega$ such that
\begin{align}
\|u\|_{H^2(\Omega)}
\leq C\|Au\|_{L^2(\Omega)}
+C\|u\|_{H^1(\Omega)}  .
\label{A1}
\end{align}
for all $u\in H^2\cap H^1_0$ such that
the support of $u$ is contained in $B_R(y)\cap\overline\Omega$.
The radius $R$ depends only on the semi-norms
$|a_{ij}|_{W^{1,N+\alpha}(\Omega)}$
and $|a_{ij}|_{C^\gamma(\overline\Omega)}$.
}
\end{lemma}
\noindent{\it Proof.}$\,\,\,$
Following the proof of Lemma 3.1.3.2 in \cite{Grisvard}
(see page 143, (3.1.3.4) and the equality above (3.1.3.5)),
we have (using our notations)
\begin{align*}
\|u\|_{H^2(\Omega)}
&\leq C\|Au\|_{L^2(\Omega)}
+C\sum_{i,j=1}^N\max_{x\in V_y}|a_{ij}(y)-a_{ij}(x)|\|u\|_{H^2(\Omega)} \\
&\quad 
+C\sum_{i,j=1}^N\|\partial_i a_{ij}\partial_j u\|_{L^2(\Omega)}\\
&\leq C\|Au\|_{L^2(\Omega)}
+CR^\beta\|u\|_{H^2(\Omega)} \\
&\quad 
+C\sum_{i,j=1}^N\|\nabla a_{ij}\|_{L^{N+\alpha}(\Omega)}
\|\nabla u\|_{L^{2(N+\alpha)/(N-2+\alpha)}(\Omega)} .
\end{align*}
When $R$ is small enough we have
\begin{align*}
\|u\|_{H^2(\Omega)}
\leq C\|Au\|_{L^2(\Omega)}
+C\|\nabla u\|_{L^{2(N+\alpha)/(N-2+\alpha)}(\Omega)} .
\end{align*}
Since $H^2$ is compactly embedded into
$W^{1,2(N+\alpha)/(N-2+\alpha)}$
which is again embedded into $H^1$,
there exists $\theta_\alpha\in (0,1)$ such that
\begin{align*}
\|\nabla u\|_{L^{2(N+\alpha)/(N-2+\alpha)}(\Omega)}
&\leq \epsilon \|u\|_{H^2(\Omega)}
+ C_\epsilon \|u\|_{H^1(\Omega)} ,\quad\forall\,\epsilon\in(0,1).
\end{align*}
Choosing $\epsilon$ small enough, \refe{A1} follows from the last
the last two inequalities
This completes the proof of Lemma \ref{LemLocH2}. \,\endproof
\medskip

Then \refe{H2Reg0} can be proved by using Lemma \ref{LemLocH2}
and a perturbation procedure (as mentioned in
the proof of \cite[Theorem 3.1.3.1]{Grisvard}).\\

\noindent{\it Proof of \refe{W1inftyReg}:}$\,\,$
Theorem 3.4 of [13] states that
if $\Omega\subset\R^N$ ($N\geq 3$) is convex
and the coefficients $a_{ij}$ are H\"older continuous (so that
(3.1)-(3.3) of  [13] hold),
the Green's function of the elliptic operator
$A$
with the Dirichlet boundary condition
satisfies
\begin{align}
| \nabla_xG_e(x,y)|+| \nabla_yG_e(x,y)|
\leq \frac{C}{|x-y|^{N-1}} ,
\label{A2}
\end{align}
where we have used the symmetry
$G_e(x,y)=G_e(y,x)$.
Therefore, any $H^1_0$ solution of the equation
$-\nabla\cdot(a\nabla u)=f$ satisfies
\begin{align*}
|\nabla u(x)|=\bigg|\int_\Omega \nabla_xG_e(x,y)f(y)\d y\bigg|
\leq \bigg\|\frac{C}{|x-y|^{N-1}}\bigg\|_{L^{p'}(\Omega)}
\|f\|_{L^p(\Omega)}
\leq C\|f\|_{L^p(\Omega)}
\end{align*}
for $p>N$. As pointed out in \cite{Fromm} (page 227,
the paragraph below Proposition 1),
the inequality \refe{A2} for $N=2$ can be proved
with some minor modifications
on the proof of \cite[Theorem 3.3--3.4]{GW}
since Theorem 3.3--3.4 of \cite{GW} only requires
$a_{ij}$ being
H\"older continuous coefficients.

\vskip0.1in

\noindent{\it Proof of \refe{LpW2quf}-\refe{LpW1quf}:}
Since $W^{1,N+\alpha}(\Omega)\hookrightarrow C(\overline\Omega)$,
Theorem 1 of \cite{JLW} implies that
the solution of the elliptic equation
\begin{align}
\left\{\begin{array}{ll}
A u = f
&\mbox{in}\,\,\,\Omega, \\
u=0 &\mbox{on}\,\,\,\partial\Omega,
\end{array}\right.
\end{align}
with continuous coefficients $a_{ij}$ in a
convex domain $\Omega\subset\R^N$ satisfies
\begin{align}\label{W1quf}
\|u\|_{W^{1,q}(\Omega)}\leq C\|f\|_{W^{-1,q}(\Omega)} ,
\quad\forall\, 1<q<\infty  ,
\end{align}
where we have noted that a continuous function
is $(\delta,R)$ vanishing, and
a convex domain is $(\delta,\sigma,R)$-quasiconvex \cite{JLW}.
Since the solution of \refe{PDE0} with $f=0$
satisfies (integrating the equation against $|u|^{q-2}u$)
$$
\|u\|_{L^q(\Omega)} \leq \|u^0\|_{L^q(\Omega)} \, ,
$$
it follows that the semigroup generated by the
elliptic operator $A$ is a
contraction semigroup on $L^q(\Omega)$.
By Theorem 1, Section 2, Chapter 3 of
\cite{Stein}, the semigroup
$\{E(t)\}_{t>0}$ generated by the elliptic operator $A$
has an analytic continuation
(analyticity of the semigroup $\{E(t)\}_{t>0}$). Moreover,
by the maximum principle
we have $u^0\geq 0$ $\implies$ $u\geq 0$
(positivity of the semigroup $\{E(t)\}_{t>0}$) and then,
by Corollary 4.d of \cite{Weis2}, the
solution of \refe{PDE0} with $u^0=0$
has the maximal $L^p$ regularity
\refe{LpW2quf}.

In other words, the map from $f$ to $Au$ given by the
formula
$$
Au=\int_0^t AE(t-s)f(\cdot,s)\d s
$$
is bounded in $L^p((0,T);L^q)$, for all $1<p,q<\infty$.
Since
$$
A^{1/2}u=\int_0^t A^{1/2}E(t-s)f(\cdot,s)\d s
=\int_0^t AE(t-s) A^{-1/2}f(\cdot,s)\d s ,
$$
it follows that
\begin{align}
\|A^{1/2}u\|_{L^p((0,T);L^{q})}
&\leq C_{p,q}\|A^{-1/2}f\|_{L^p((0,T);L^{q})}  ,\quad\forall\,
1<p,q<\infty .
\end{align}
It remains to prove the boundedness of the Riesz transform
$\nabla A^{-1/2}$:
\begin{align}\label{RiezTrf}
\|\nabla A^{-1/2}f\|_{L^{q}}
&\leq C_{q}\|f\|_{L^{q}}  ,\quad\forall\,
1<q<\infty .
\end{align}
Then the last two inequalities imply
\begin{align*}
\|\nabla u\|_{L^p((0,T);L^{q})}
&=\|\nabla A^{-1/2}(A^{1/2}u)\|_{L^p((0,T);L^{q})} \\
&\leq C_q\|A^{1/2}u\|_{L^p((0,T);L^{q})} \\
&\leq C_{p,q}\|A^{-1/2}f\|_{L^p((0,T);L^{q})}  \\
&\leq C_{p,q}\|f\|_{W^{-1,q}} ,
&&\forall\, 1<p,q<\infty ,
\end{align*}
where the last step of the inequality above is due to the following
duality argument ($A^{-1/2}$ is self-adjoint):
\begin{align*}
(A^{-1/2}f,g)  =(f, A^{-1/2}g)
\leq C\|f\|_{W^{-1,q}}\|\nabla A^{-1/2}g\|_{L^{q'}}
\leq C\|f\|_{W^{-1,q}}\|g\|_{L^{q'}} .
\end{align*}

It has been proved in \cite[Theorem B]{Shen05}
that the Riesz transform is bounded on $L^q(\Omega)$
(i.e. the inequality (\ref{RiezTrf}) holds) if and only if the solution of
the homogeneous equation
\begin{align}\label{HomoEqu}
A u = 0
\end{align}
satisfies the local estimate 
\begin{align}\label{HomoEst0}
\bigg(\frac{1}{r^N}\int_{\Omega\cap B_r(x_0)}
|\nabla u|^q\d x\bigg)^{\frac{1}{q}}
\leq C
\bigg(\frac{1}{r^N}\int_{\Omega\cap B_{\sigma_0r}(x_0)}
|\nabla u|^2\d x\bigg)^{\frac{1}{2}} 
\end{align}
for all $x_0\in\Omega$ and $0<r<r_0$, 
where $r_0$ and $\sigma_0\geq 2$
are any given small positive constants
such that $\Omega\cap B_{\sigma_0r_0}(x_0)$ can be given
by the intersection of $B_{\sigma_0r_0}(x_0)$ with a Lipschitz graph.
It remains to prove \refe{HomoEst0}.

Let $\omega$ be a smooth cut-off function which equals
zero outside $B_{2r}:=B_{2r}(x_0)$
and equals 1 on $B_r$.
Extend $u$ to be zero on $B_{2r}\backslash\Omega$
and denote by $u_{2r}$ the average of $u$ over $B_{2r}$.
Then (\ref{HomoEqu}) implies
\begin{align}\label{HomoEqu2}
&\sum_{i,j=1}^N\partial_i(a_{ij}\partial_j (\omega (u-u_{2r}))) \nn \\
&=\sum_{i,j=1}^N\partial_i(a_{ij}(u-u_{2r})\partial_j\omega )
+\sum_{i,j=1}^N a_{ij}\partial_i\omega\partial_j (u-u_{2r})
\quad\mbox{in}\,\,\,\Omega,
\end{align}
and the $W^{1,q}$ estimate (\ref{W1quf})
implies
\begin{align*}
\|\omega (u-u_{2r})\|_{W^{1,q}(\Omega)}
&\leq C\|(u-u_{2r})\partial_j\omega\|_{L^q(\Omega)}
+C\|\partial_i\omega\partial_j u\|_{W^{-1,q}(\Omega)} \\
&\leq C\|(u-u_{2r})\partial_j\omega\|_{L^q(\Omega)}
+C\|\partial_i\omega\partial_j u\|_{L^s(\Omega)} \\
&= C\|(u-u_{2r})\partial_j\omega\|_{L^q(B_{2r})}
+C\|\partial_i\omega\partial_j u\|_{L^s(B_{2r})}  \\
&\leq Cr^{-1}\|\nabla u\|_{L^s(B_{2r})} ,
\end{align*}
where $s=qN/(q+N)<q$ satisfies
$L^s(\Omega)\hookrightarrow W^{-1,q}(\Omega)$
and $W^{1,s}(\Omega)\hookrightarrow L^q(\Omega)$.
The last inequality implies
\begin{align}\label{HomoEst2}
\|\nabla u\|_{L^q(\Omega\cap B_{r})}
&\leq Cr^{-1}\|\nabla u\|_{L^s((\Omega\cap B_{2r})} .
\end{align}
If $s\leq 2$ then one can derive
\begin{align*}
\|\nabla u\|_{L^q(\Omega\cap B_{r})}
&\leq Cr^{N/q-N/2}\|\nabla u\|_{L^2(\Omega\cap B_{2r})} .
\end{align*}
by using one more H\"older's inequality on the right-hand side.
Otherwise,
one only needs a finite number of iterations of (\ref{HomoEst2})
to reduce $s$ to be less than $2$.
This completes the proof of (\ref{HomoEst0}).

The proof of \refe{LpW2quf}-\refe{LpW1quf} is complete.
\,\endproof

\bigskip

\noindent{\bf Acknowledgement}~~
We would like to thank the anonymous referees
for many valuable comments and suggestions, 
which are very helpful to improve
both the quality and presentation of this paper.

\end{document}